\documentclass[a4paper,11pt]{article}
\usepackage[utf8]{inputenc}
\usepackage[T1]{fontenc}
\usepackage{lmodern}

\usepackage{boxedminipage}
\usepackage{amsfonts}
\usepackage{amsmath} 
\usepackage{amssymb}
\usepackage{graphicx}
\usepackage{amsthm}
\usepackage{subfig}
\usepackage{cancel}
\usepackage{microtype}
\usepackage{hyperref}
\newtheorem{theorem}{Theorem}[section]
\newtheorem{lemma}[theorem]{Lemma}

\newtheorem{corollary}[theorem]{Corollary}
\newtheorem{definition}[theorem]{Definition}

\newtheorem{fact}[theorem]{Fact}
\newtheorem{question}{Question}

\newtheorem*{theorem*}{Theorem}

\newcommand{\forceP}{\mathbb{P}}
\newcommand{\forceQ}{\mathbb{Q}}

\newcommand{\ZFC}{\mathsf{ZFC}}

\newcommand{\ZFCminus}{\mathsf{ZFC}^{-}}

\newcommand{\CH}{\mathsf{CH}}

\newcommand{\PD}{\mathsf{PD}}

\def\undertilde#1{\mathord{\vtop{\ialign{##\crcr
$\hfil\displaystyle{#1}\hfil$\crcr\noalign{\kern1.5pt\nointerlineskip}
$\hfil\tilde{}\hfil$\crcr\noalign{\kern1.5pt}}}}}

\title{Forcing the ${\Pi^1_n}$-Uniformization Property}
\author{Stefan Hoffelner\footnote{The author's research was funded in part by the Austrian Science Fund (FWF) Grant-DOI 10.55776/P37228. This research was also funded in part by the Deutsche Forschungsgemeinschaft (DFG German Research Foundation) under Germany's Excellence Strategy EXC 2044 390685587, Mathematics M\"unster: Dynamics--Geometry--Structure. For the purpose of open access, the author has applied a CC BY public copyright license to any Author Accepted Manuscript version arising from this submission. } }
\date{November 9, 2026}
\begin{document}

\maketitle

\begin{abstract}
We generically construct a model in which the ${\Pi^1_3}$-uniformization property is true, thus lowering the best known consistency strength from the existence of $M_1^{\#}$ to just $\ZFC$. The forcing construction also works over canonical inner models with Woodin cardinals, which yields, for the first time, universes where the $\Pi^1_{2n}$-uniformization property holds for $n>1$, thus producing models which contradict the natural $\PD$-induced pattern.

\end{abstract}

\section{Introduction}
The question of finding nicely definable 
choice functions for a definable family of sets is an old and well-studied subject in descriptive set theory; see, for instance, \cite{Kechris,Moschovakis}. The uniformization problem, first mentioned by N. Lusin in 1930 (see \cite{Lusin}), asks to find choice functions which lie at the same projective level as the set they aim to uniformize.
Recall that for $A\subseteq 2^{\omega}\times 2^{\omega}$, a
\emph{uniformization} of $A$ is a function
\[
   f:\operatorname{pr}_1(A)\longrightarrow 2^\omega
\]
such that
\[
   \operatorname{graph}(f)\subseteq A.
\]
Equivalently, for every $x\in\operatorname{pr}_1(A)$ one has
$(x,f(x))\in A$.

\begin{definition}
 We say that a pointclass $\Gamma$ has the uniformization
 property iff every element of $\Gamma$ admits a uniformization 
 in $\Gamma$.
\end{definition}

It is a classical result due to M. Kondo \cite{Kondo} that lightface $\Pi^1_{1}$-sets do have the 
uniformization property, this also yields the uniformization property for $\Sigma^1_{2}$-sets. This is as much as $\ZFC$ can prove about uniformization.
In the constructible universe $L$, as shown by J. Addison in \cite{Addison}, $\Sigma^1_{n}$ does have the uniformization property for $n \ge 3$, which follows from the existence of a $\Sigma^1_{2}$-good wellorder of the reals; for background on projective wellorders see also \cite{CS,SyVera}. Thus the $\Pi^1_{n}$-uniformization fails for $n \ge 3$. On the other hand, by the celebrated results of Y. Moschovakis (see \cite{Moschovakis2}, Theorem 1), 
$\bf{\Delta}^1_{2n}$-projective determinacy implies ${\Pi}^1_{2n+1}$-uniformization, yet the determinacy assumption exceeds in logical strength $\ZFC$. It is known due to W. H. Woodin (see \cite{MSW}), that $\bf{\Delta}^1_2$-projective determinacy implies the existence of $M_1^{\#}$, hence yields an inner model with a Woodin cardinal. 
As with other regularity properties of the reals like Lebesgue measurability or Baire property, which both follow from $\PD$ as well, it is natural to ask whether the $\Pi^1_3$-uniformization property bears large cardinal strength as well. We shall answer it negatively.

\begin{theorem*}
There is a generic extension of $L$ in which the $\Pi^1_3$-uniformization property is true. 
\end{theorem*} 
The same uniformization construction also applies over canonical inner models with Woodin cardinals. This can be used to obtain better upper bounds in terms of consistency strength for the $\Pi^1_n$ uniformization property for odd $n$. For even $n\geq 4$ we can produce for the first time models where the $\Pi^1_n$ uniformization property holds true.
\begin{theorem*}
Let $M_n$ be the canonical inner model with $n$ Woodin cardinals. Then there is a generic extension of $M_n$ in which the $\Pi^1_{n+3}$ uniformization property holds true.
\end{theorem*}

Since the first version of this paper was written, the method developed here has led to further constructions.  One combines global $\Sigma$-uniformization with a large continuum and a projective well-order \cite{HoffelnerLargeContinuumGlobalSigma}, and the $\Sigma$-side of the local Delfino problem combines related coding methods with projective regularity and uniformization \cite{HoffelnerLocalDelfinoSigma}.  These later papers refine and apply the forcing method introduced here.  The present paper contains the fixed point construction and the hierarchy of allowable iterations needed to force projective $\Pi$-uniformization without assuming the corresponding determinacy hypothesis.

Questions concerning the forceability of (local) consequences of $\PD$ do have a long tradition in set theory. There is a vast body of literature concerning the forcability of local levels of the projective hierarchy satisfying (Boolean combinations of) the Baire property, the perfect set property or Lebesgue measurability; see for instance \cite{FS}. There has been very little progress in the past, however, concerning similar questions for the separation, the reduction and the uniformization property. Indeed, even the question of whether one can force the $\Sigma^1_3$-separation property, which is the weakest of said properties, remained an open problem for 50 years (see \cite{Mathias}, Problem 3029).

This article continues this line of research and develops the forcing construction used in the subsequent work mentioned above. It is organized as follows: in the preliminaries section, we briefly introduce the forcings which we will use in the proof and produce a generic extension $W$ of $L$ which will be a well-suited ground model for our needs. We then start to prove the theorems from above. The main idea is to turn the problem of finding a partial order $\forceP$ which forces the $\Pi^1_3$-property into a fixed point problem.  We shall define a derivation operator which acts on a specific set of forcing iterations of length $\omega_1$. This operator will be applied transfinitely often, and will produce better and better approximations to the set of forcings we actually want to use in the end. The process is shown to converge in that eventually a fixed point $\mathcal{P}$, i.e. a suitable set of forcings is reached. Forcings which belong to this fixed point $\mathcal{P}$ allow for a certain, seemingly self-referential  definition of an $\omega_1$-length iteration which forces the $\Pi^1_3$-uniformization property over $W$.
We then follow up, to alter the said process such that it becomes applicable to the  canonical inner models $M_n$ with $n$ Woodin cardinals.

There are some similarities to \cite{H}, in particular the two proofs rely on a similar ground model $W$, which is a generic extension of $L$, and use a similar coding method which relies on a suitably chosen $\omega_1$-sequence of $\omega_1$-Suslin trees. However a more straightforward application of the ideas of \cite{H} will fail badly to produce a model of the $\Pi^1_3$-uniformization property. As a consequence, a solution has to necessarily introduce several new ideas in order to succeed. The presentation of those is the goal of this paper.
\section{Preliminaries}

\subsection{Notation}
The notation we use will be mostly standard, we hope. We write $\forceP=(\forceP_{\alpha} \, : \, \alpha < \gamma)$ for a forcing iteration of length $\gamma$ with initial segments $\forceP_{\alpha}$. The $\alpha$-th factor of the iteration will be denoted with $\forceP(\alpha)$. Note here that we drop the dot on $\forceP(\alpha)$, even though $\forceP(\alpha)$ is in fact a $\forceP_{\alpha}$-name of a partial order.

\subsection{The forcings which are used}
The forcing notions used in the construction are standard. The preparatory extension is
obtained solely by adding branches through a fixed sequence of Suslin trees,
and the later coding uses the ground-model version of the $\omega_1$-Cohen
forcing together with Jensen--Solovay almost disjoint coding \cite{JensenSolovay}.

Recall that a Suslin tree is a normal tree of height $\omega_1$ with countable
levels and without an uncountable chain or antichain.  Whenever $S$ is a Suslin
tree, we write $\mathbb{S}$ for the forcing whose conditions are the nodes of $S$,
ordered by end-extension.  Thus forcing with $\mathbb{S}$ adds a cofinal branch through
$S$.

\begin{definition}
Let $\vec{S}=\langle S_\xi\mid \xi<\kappa\rangle$ be a sequence of Suslin
trees.  We call $\vec{S}$ \emph{independent} if, for every nonempty finite set
$e\subseteq\kappa$, the product
\[
   \prod_{\xi\in e}S_\xi
\]
is again a Suslin tree.
\end{definition}

For $A\subseteq\kappa$ let
\[
   \operatorname{Br}_A(\vec S)=\prod_{\xi\in A}^{\mathrm{fs}}\mathbb{S}_\xi
\]
denote the finite-support product of the corresponding branch forcings.  When
the sequence is fixed, we simply write $\operatorname{Br}_A$.

\begin{lemma}\label{finite-support-branch-products}
Let $\vec S=\langle S_\xi\mid\xi<\kappa\rangle$ be an independent sequence of
Suslin trees, where $\kappa\leq\omega_1$.
\begin{enumerate}
\item For every $A\subseteq\kappa$, the forcing $\operatorname{Br}_A$ is c.c.c.
\item If $A\subseteq B\subseteq\kappa$, then $\operatorname{Br}_A$ is a
      complete subforcing of $\operatorname{Br}_B$, and
      \[
         \operatorname{Br}_B\cong
         \operatorname{Br}_A\times\operatorname{Br}_{B\setminus A}.
      \]
\item If $\xi\notin A$, then
      \[
         \operatorname{Br}_A\Vdash ``S_\xi\text{ is a Suslin tree}."
      \]
\end{enumerate}
\end{lemma}

\begin{proof}
For (1), let $\{p_\alpha\mid\alpha<\omega_1\}$ be a family of conditions in
$\operatorname{Br}_A$.  After thinning it out, the supports form a
$\Delta$-system with finite root $r$, and the restrictions
$p_\alpha\upharpoonright r$ belong to the finite product
$\prod_{\xi\in r}S_\xi$.  By independence this finite product is Suslin, hence
c.c.c., so two of the root restrictions are compatible.  Since the remaining
parts of their supports are disjoint, the corresponding two conditions in
$\operatorname{Br}_A$ are compatible.

The factorization in (2) is the usual coordinatewise factorization of a
finite-support product.  For (3), observe that
\[
   \operatorname{Br}_A\times\mathbb{S}_\xi
   \cong \operatorname{Br}_{A\cup\{\xi\}},
\]
which is c.c.c. by (1).  The standard product characterization of preservation
of a Suslin tree now yields the assertion.
\end{proof}

We next recall the almost disjoint coding forcing due to Jensen and Solovay \cite{JensenSolovay}.
We identify subsets of $\omega$ with their characteristic functions.  Let
$D=\langle d_\alpha\mid\alpha<\omega_1\rangle$ be an almost disjoint family of
infinite subsets of $\omega$, and let $X\subseteq\omega_1$.

\begin{definition}\label{definitionadcoding}
The almost disjoint coding forcing $\mathbb A_D(X)$ consists of pairs
$(r,R)\in[\omega]^{<\omega}\times[D]^{<\omega}$.  We put
$(s,S)\leq(r,R)$ if and only if
\begin{enumerate}
\item $r\subseteq s$ and $R\subseteq S$;
\item whenever $\alpha\in X$ and $d_\alpha\in R$, then
      $s\cap d_\alpha=r\cap d_\alpha$.
\end{enumerate}
\end{definition}

If $G\subseteq\mathbb A_D(X)$ is generic and
$x_G=\bigcup\{r:(r,R)\in G\}$, then
\[
   \alpha\in X\quad\Longleftrightarrow\quad
   x_G\cap d_\alpha\text{ is finite}.
\]
Fix a recursive injection $\#:2^{<\omega}\to\omega$. If
$\langle x_\xi:\xi<\omega_1^L\rangle$ enumerates the constructible
reals in the canonical order, use
$d_\xi=\{\#(x_\xi\restriction k):k<\omega\}$ and
$D=\langle d_\xi:\xi<\omega_1^L\rangle$.

\begin{lemma}\label{a.d.coding preserves Suslin trees}
For every $X\subseteq\omega_1$, the forcing $\mathbb A_D(X)$ is Knaster.  In
particular it is c.c.c. and preserves every ground-model Suslin tree.
\end{lemma}

\begin{proof}
The Knaster property follows from the usual $\Delta$-system argument applied
to the finite side conditions.  If $T$ is a Suslin tree, then the product of a
Knaster forcing with the c.c.c. forcing $\mathbb{T}$ is c.c.c.  Hence
$\mathbb A_D(X)\times\mathbb{T}$ is c.c.c., which is equivalent to preservation of
$T$ as a Suslin tree.
\end{proof}

Finally, we write
\[
   \mathbb C_L(\omega_1)=({}^{<\omega_1}2)^L
\]
for the $L$-computed $\omega_1$-Cohen forcing, ordered by end-extension.  Thus
conditions are countable binary sequences belonging to $L$.  Whenever this
forcing is used in an outer model of $L$, the notation refers to this fixed
$L$-poset rather than to a recomputation in the ambient universe.

\subsection[The ground model W of the iteration]{The ground model $W$ of the iteration}
We now define the ground model over which the main iteration will be carried
out.  We start with the canonical independent Suslin sequence already present
in $L$ and add a cofinal generic branch through every member of that sequence.

The following lemma gives the simultaneous form of Jensen's diamond
construction which we use throughout the paper; see also \cite{Jech}.  We use
the levelwise product of finitely many trees.  This is a dense suborder of the
corresponding finite forcing product.

\begin{lemma}
\label{canonical-L-Suslin-sequence}
There is a canonical $L$-definable sequence
\[
   \vec S=\langle S_\xi\mid\xi<\omega_1\rangle
\]
of independent normal Suslin trees.  The construction is given by a recursion
which is uniform over $L_{\omega_1}$.
\end{lemma}

\begin{proof}
Work in $L$.  Fix the canonical $L$-diamond sequence
$\langle a_\alpha\mid\alpha<\omega_1\rangle$ and the canonical $L$-coding of
finite sets of ordinals, countable trees, and subsets of their finite products.
We recursively construct coherent countable approximations
$\vec S^{\,\alpha}=\langle S^\alpha_\xi\mid \xi<\alpha\rangle$ such that each
$S^\alpha_\xi$ is a normal splitting tree of height $\alpha$ and
$S^\beta_\xi=S^\alpha_\xi\restriction\beta$ whenever $\xi<\beta<\alpha$.  At
limit stages we take unions.

It remains to describe the passage from $\alpha$ to $\alpha+1$.  For old
coordinates $\xi<\alpha$, successor stages are handled by putting countably
many immediate successors above each node on the top level.  Suppose now that
$\alpha$ is a nonzero limit.  If $a_\alpha$ does not decode as a pair
$(e,A)$, where $e\in[\alpha]^{<\omega}\setminus\{\varnothing\}$ and $A$ is a
maximal antichain in the product
\[
   \mathbb{P}_e^\alpha=\bigcup_{\delta<\alpha}
      \prod_{\xi\in e}(S^\alpha_\xi)_\delta,
\]
then we use the $<_L$-least countable family of cofinal branches through each
old tree which covers all old nodes, and put one new top node above each of
these branches.

If $a_\alpha$ does decode such a pair $(e,A)$, we choose the branch families
for the coordinates in $e$ simultaneously so as to seal $A$.  Concretely, by a
countable genericity argument in
$\prod_{(\xi,n)\in e\times\omega}^{\mathrm{fs}}\mathbb{S}^\alpha_\xi$, choose cofinal
branches $b_{\xi,n}$, $\xi\in e$, $n<\omega$, such that they cover all old
nodes and such that, for every $\nu:e\to\omega$, the tuple of new top nodes
placed above $\langle b_{\xi,\nu(\xi)}\mid\xi\in e\rangle$ extends some member
of $A$.  The required dense sets are just the usual ones ensuring cofinality,
coverage of old nodes, and the displayed sealing requirement; density of the
last kind uses maximality of $A$ in $\mathbb{P}_e^\alpha$.  For $\xi\notin e$ we again
use the default $<_L$-least covering family.  Finally, the new coordinate
$\alpha$ is initialized as the canonical finite-support binary tree of height
$\alpha+1$.  This completes the recursion.  Put
$S_\xi=\bigcup_{\xi<\alpha<\omega_1}S^\alpha_\xi$.  Then each $S_\xi$ is a
normal splitting tree of height $\omega_1$ with countable levels.

We show that the sequence is independent.  Fix a finite nonempty
$e\subseteq\omega_1$ and let
\[
   \mathbb{P}_e=\bigcup_{\alpha<\omega_1}
      \prod_{\xi\in e}(S_\xi)_\alpha
\]
be the product.  This is dense in the usual finite product of the
corresponding tree forcings.  Let $A\subseteq \mathbb{P}_e$ be a maximal antichain in
$L$.  There is a club of $\alpha>\sup(e)$ such that
$A\cap(\mathbb{P}_e\restriction\alpha)$ is maximal in $\mathbb{P}_e\restriction\alpha$, and,
by coherence of the coding, a club on which the fixed code $x_{e,A}$ for
$(e,A)$ restricts to a code for
$(e,A\cap(\mathbb{P}_e\restriction\alpha))$.  By $\diamondsuit$, choose a limit
$\alpha$ in the intersection of these clubs with
$a_\alpha=x_{e,A}\cap\alpha$.  At stage $\alpha$ the construction seals
$A\cap(\mathbb{P}_e\restriction\alpha)$.  Hence every node on level $\alpha$, and
therefore every later node, extends a member of this initial antichain.  Since
$A$ itself is an antichain, $A=A\cap(\mathbb{P}_e\restriction\alpha)$, so $A$ is
countable.  Thus $\mathbb{P}_e$ has no uncountable antichain.  A cofinal branch through
$\mathbb{P}_e$ would, by normal splitting, yield an uncountable antichain by choosing at
successor levels nodes incompatible with the branch.  Hence $\mathbb{P}_e$ is Suslin.
Since $e$ was arbitrary, $\vec S$ is independent.  The construction is
canonical from $<_L$, the diamond sequence, and the fixed coding, and is
therefore uniform over $L_{\omega_1}$.
\end{proof}

\begin{lemma}
\label{definabilityofvecS}
Let $M$ be a transitive uncountable model of
\[
   \ZFCminus+``\aleph_1\text{ exists}"
\]
such that
\[
   \omega_1^M=\omega_1^L=\omega_1.
\]
There is a $\Sigma_1$-formula
\[
   \vartheta_{\mathrm{Sus}}(t,\gamma,\eta,\omega_1),
\]
independent of $M$, which defines over $M$ the relation
\[
   \{(t,\gamma,\eta):
      \gamma,\eta<\omega_1\text{ and }t\in(S_\eta)_\gamma\}.
\]
Equivalently, there is a formula
$\operatorname{CanSuslin}_L(X,\omega_1)$, which may be taken to be
$\Sigma_1(\omega_1)$ over $H(\omega_2)$, such that in every
$\omega_1$-preserving outer model $U$ of $L$,
\[
   U\models\operatorname{CanSuslin}_L(X,\omega_1)
   \quad\Longleftrightarrow\quad X=\vec S.
\]
In particular, the canonical sequence $\vec S$ is computed in the same way in
all the outer models considered below.
\end{lemma}

\begin{proof}
By Lemma~\ref{canonical-L-Suslin-sequence}, the sequence \(\vec S\) is
obtained by a canonical recursion over \(L_{\omega_1}\), and the corresponding
level relation is \(\Sigma_1(\omega_1)\)-definable over \(L_{\omega_1}\).
If \(M\) is transitive, satisfies \(\ZFCminus+``\aleph_1\text{ exists}"\),
and has \(\omega_1^M=\omega_1^L\), then \(M\) computes \(L_{\omega_1}\)
correctly.  Hence \(M\) computes the same canonical diamond sequence, the same
coding apparatus, and therefore the same recursive construction of
\(\vec S\).  The same \(\Sigma_1\)-formula
\(\vartheta_{\mathrm{Sus}}(t,\gamma,\eta,\omega_1)\) consequently defines over
\(M\) the relation \(t\in(S_\eta)_\gamma\).

The formulation by \(\operatorname{CanSuslin}_L(X,\omega_1)\) is obtained by
saying that \(X\) is the union of the unique system produced by this
canonical recursion.  Since every \(\omega_1\)-preserving outer model of \(L\)
has the same \(L_{\omega_1}\), it satisfies this formula exactly for the true
sequence \(\vec S\).
\end{proof}

Let
\[
   \operatorname{Br}=\prod_{\xi<\omega_1}^{\mathrm{fs}}\mathbb{S}_\xi
\]
be the finite-support branch product, computed in $L$, and let
$H\subseteq\operatorname{Br}$ be $L$-generic.  For every $\xi<\omega_1$ set
\[
   b_\xi=\bigcup\{p(\xi):p\in H\text{ and }\xi\in\operatorname{dom}(p)\}.
\]
Then $b_\xi$ is a cofinal $L$-generic branch through $S_\xi$.  In this paper,
to \emph{destroy} a Suslin tree always means to add such a cofinal generic
branch through it.  We define
\[
   W=L[H].
\]

\begin{theorem}\label{properties-of-W}
The model $W$ has the following properties.
\begin{enumerate}
\item $\operatorname{Br}$ is a c.c.c. forcing of size $\aleph_1$.  Hence $W$
      preserves all cardinals and cofinalities of $L$, and $W\models\CH$.
\item Every $S_\xi$ is destroyed in $W$, witnessed by the branch $b_\xi$.
\item For $A\subseteq\omega_1$ in $L$, let
      $H_A=H\cap\operatorname{Br}_A$.  Then $H_A$ is
      $\operatorname{Br}_A$-generic over $L$, and
      \[
         L[H_A]\models
         ``S_\xi\text{ has a cofinal branch}"
         \quad\Longleftrightarrow\quad \xi\in A.
      \]
\item Every real in $W$ belongs to $L[H_A]$ for some countable
      $A\subseteq\omega_1$ belonging to $L$.
\end{enumerate}
\end{theorem}

\begin{proof}
The first assertion follows from Lemma~\ref{finite-support-branch-products}.
Since $L\models\CH$ and $\operatorname{Br}$ has size $\aleph_1$, there are only
$\aleph_1$ many nice $\operatorname{Br}$-names for reals, and therefore $W$
still satisfies $\CH$.  The definition of $H$ and the dense sets requiring a
condition to reach a prescribed level of a prescribed coordinate give (2).

For (3), use the factorization
\[
   \operatorname{Br}\cong
   \operatorname{Br}_A\times\operatorname{Br}_{\omega_1\setminus A}.
\]
If $\xi\in A$, then $b_\xi\in L[H_A]$.  If $\xi\notin A$, then
Lemma~\ref{finite-support-branch-products} shows that
$\operatorname{Br}_A$ preserves $S_\xi$ as a Suslin tree, so no cofinal branch
through $S_\xi$ belongs to $L[H_A]$.

Finally, let $\dot x$ be a nice $\operatorname{Br}$-name for a real.  The
antichains occurring in $\dot x$ are countable, and the union of the finite
supports of all conditions appearing in them is therefore a countable set
$A\subseteq\omega_1$ in $L$.  Thus $\dot x$ is already a
$\operatorname{Br}_A$-name, proving (4).
\end{proof}

The relevant inner models of $W$ used below are obtained by retaining
selected branch coordinates from $H$.

\paragraph{Coding forcing.}
The extensions under consideration preserve $\omega_1^L$, so
$\omega_1=\omega_1^L$. Let $u\in 2^\omega$ and let $g$ be the generic
subset of $\omega_1$ added at a new $\mathbb C_L(\omega_1)$-coordinate.
Let $\varrho_L:({}^{<\omega_1}2)^L\to\omega_1$ be the canonical bijection. Put
\[
 \begin{aligned}
 h(g)&=\{\varrho_L(g\restriction\gamma):\gamma<\omega_1\},\\
 I(u,g)&=\{\omega\gamma+2k+u(k):\gamma\in h(g),\ k<\omega\}.
 \end{aligned}
\]
Let $X\subseteq\omega_1$ canonically code
$(u,h(g),\langle b_\xi:\xi\in I(u,g)\rangle)$. We describe the
localization $Y=Y(X)$ used by $\operatorname{Code}(u,\eta)$,
where $\eta$ specifies the coordinate supplying $g$.

Use a recursive theory $T$ consisting of $\ZFCminus$ with Collection
and Choice, existence of $\omega_1$, the relativizations of these
axioms to $L[a]$ for each real $a$, and
$\operatorname{HC}^L=L_{\omega_1^L}$. For a transitive $N\models T$
with $u\in N$ and $\omega_1^N=(\omega_1^L)^N=\nu$, put
\[
 \operatorname{Pat}(N,u,\nu)
 \quad\Longleftrightarrow\quad
 N\models \forall\beta<\nu\ \exists\gamma\in[\beta,\nu)\ \operatorname{Pattern}(\gamma,u).
\]
Here $\operatorname{Pattern}(\gamma,u)$ asserts that for every
$k<\omega$ there is a cofinal branch through $S_{\omega\gamma+2k+u(k)}$. The predicate includes the displayed requirements
on $N,u,\nu$. The trees are defined by Lemma~\ref{definabilityofvecS}.

In $L[X]$, collapse an elementary submodel of a sufficiently large
$H_\chi$ of size $\aleph_1$ containing $X,u,\omega_1$ and all ordinals
below $\omega_1$. Its collapse $N$ satisfies
$\operatorname{Pat}(N,u,\omega_1)$. Choose a bijection $e:\omega_1\to N$
and code membership and satisfaction via $e$ by $Z\subseteq\omega_1$,
marking the indices of $u$ and $\omega_1$. Use finite tuple codes ordered
by maximum. In $L[Z]$, take the club of closure points $\nu$ at which $Z\cap\nu$ codes membership and satisfaction of $e[\nu]$ via $e\restriction\nu$, $e[\nu]\prec N$, with $u,\omega_1\in
e[\nu]$ and $e[\nu]\cap\omega_1=\nu$. Such a club is obtained by
closing under the Skolem maps, $e$ and its inverse. At these points the
collapse $N_\nu$ satisfies $\operatorname{Pat}(N_\nu,u,\nu)$. Let
$E$ consist of the nonzero limit points of this club and $0$. If $c\in
E$, write $d_c=\min(E\setminus(c+1))$ and let $w_c$ be the $<_L$-least
real coding a well-order of type $d_c$. Set
\[
 Y(X)=\{2\cdot\xi+1:\xi\in Z\}
       \cup\{2\cdot(c+k):c\in E,\ k<\omega,\ w_c(k)=1\}.
\]
The intervals $[c,c+\omega)$ are disjoint. Make the choices in
$L[X]$ and $L[Z]$ by their canonical well-orders, so the same $X$ gives the same $Y$ in all models $M$ with $\omega_1^M=\omega_1$. This is the
localization argument of David~\cite{David}. The second factor of
$\operatorname{Code}(u,\eta)$ is $\mathbb A_D(Y(X))$.

Let $\Psi(r,u)$ assert that for every countable transitive $M\models T$
with $r,u\in M$ and $\omega_1^M=(\omega_1^L)^M$, one has
\[
 M\models\exists N\in L[r]\ \operatorname{Pat}(N,u,\omega_1^M).
\]
A real coding $M$ can also code its satisfaction relation. The
well-foundedness requirement is $\Pi^1_1$ and the remaining tests
on the code are arithmetic. Hence $\Psi$ is $\Pi^1_2$, and
$\Phi(u)\iff\exists r\,\Psi(r,u)$ is $\Sigma^1_3$.

\begin{lemma}\label{lem:L-decoder}
The generic real $r$ of the almost disjoint factor of
$\operatorname{Code}(u,\eta)$ satisfies $\Psi(r,u)$. Conversely, in
an extension preserving $\omega_1$, if $\Psi(r,u)$ holds, then for
unboundedly many $\gamma<\omega_1$ and every $k<\omega$, $L[r]$
contains a cofinal branch through $S_{\omega\gamma+2k+u(k)}$.
\end{lemma}
\begin{proof}
Given a countable transitive $M$ as in the definition of $\Psi$,
put $\nu=\omega_1^M$. The definition of $D$ and
$\operatorname{HC}^L=L_{\omega_1^L}$ give the correct $D\restriction\nu$
in $M$. The almost disjoint real decodes $Y\cap\nu$ in $L[r]^M$.
If $\nu\notin\operatorname{Lim}(E)$, let $c=\max(E\cap\nu)$. Then
$c+\omega<\nu\le d_c$. The even part gives $w_c\in M$, contrary to
$\nu=\omega_1^M$. At a limit point of $E$ the odd part yields the code
of $N_\nu$, so the collapse lies in $L[r]^M$ and is as required.

For the converse, use $F=H_\chi^{L[r,u]}$ for sufficiently large regular
$\chi$. This model satisfies $T$ and has $\omega_1^F=(\omega_1^L)^F=\omega_1$.
If the conclusion in the definition of $\Psi$ failed at $F$,
a countable elementary submodel of $F$ containing
$r,u,\omega_1$ would collapse to a countable transitive model witnessing failure.
Thus $F$ contains such an $N\in L[r]$. The trees at $\omega_1$ are computed correctly,
and the branches in $N$ belong to $L[r]$. Finally, once $r$ witnesses
$\Phi(u)$ it continues to witness it by Shoenfield absoluteness.
\end{proof}

\subsection{Coding iterations and their initial segments}
\label{subsec:coding-iterations}

For $E\in L$ of size at most $\aleph_1$, put
\[
 \mathbb{C}_E=\operatorname{Fn}_{<\omega_1}(E\times\omega_1,2)^L.
\]
Use the Cohen index sets $E_i=[\omega_1\cdot i,\omega_1\cdot(i+1))$ and
put $E_{<\lambda}=\bigcup_{i<\lambda}E_i$. At stage $i$ we first
force with $\mathbb{C}_{E_i}$ and then with an almost disjoint factor
$\dot {\mathbb{A}}_i$, which may be trivial. Thus $\mathbb{P}_i$ is the initial segment
through the first $i$ stages, and $\mathbb{P}_i*\mathbb{C}_{E_i}$ is the forcing just
before $\dot {\mathbb{A}}_i$.

Whether to code at stage $i$, the real $u_i$ to be coded, the index
$\eta_i\in E_i$, and the isomorphism
$\varphi_i:\mathbb{C}_{E_i}\cong \mathbb{C}_{E_i}\restriction c_i$ are given by
$\mathbb{P}_i$-names. Each possible value of $\varphi_i$ is in $L$ and is
induced by bijections on the individual Cohen coordinates, as in
the proof below, with at most countably many nonidentity coordinates.
If $K_i$ is $\mathbb{C}_{E_i}$-generic, let $g_i$ be the $\eta_i$-th Cohen
subset of $\omega_1$ obtained from the filter generated by
$\varphi_i[K_i]$. At a coding stage put
\[
 I_i=I(u_i,g_i),\qquad
 X_i=\operatorname{code}
     \bigl(u_i,h(g_i),\langle b_\xi:\xi\in I_i\rangle\bigr)
\]
and use $\mathbb A_D(Y(X_i))$, possibly restricted below a condition
$a_i$. The name for $a_i$ is a $\mathbb{P}_i*\mathbb{C}_{E_i}$-name. At other stages
we still use $\mathbb{C}_{E_i}$, but $\dot {\mathbb{A}}_i$ is trivial and $I_i=\varnothing$.

For $\lambda\leq\omega_1$, write the forcing as
\begin{equation}\label{eq:coding-iteration}
 \mathbb{P}_\lambda=\mathbb{C}_{E_{<\lambda}}*\dot {\mathbb{R}}_\lambda,
 \qquad \mathbb{R}_\lambda=\mathop{*}_{i<\lambda}^{\mathrm{fs}}\dot {\mathbb{A}}_i .
\end{equation}
Here the names for $\dot {\mathbb{A}}_i$ are translated from $\mathbb{P}_i*\mathbb{C}_{E_i}$.
The later Cohen forcings are ground-model partial orders and the
names for the earlier factors use only the earlier Cohen coordinates.
Moreover the order on $\mathbb A_D(Y)\restriction a$ is determined by
$Y$ and $a$. Thus the later Cohen forcings can be moved before the
earlier almost disjoint factors. At a limit,
\eqref{eq:coding-iteration} uses the countable-support Cohen product
computed in $L$ and the finite-support almost disjoint iteration.
For $\lambda=\omega_1$, the Cohen product is isomorphic to $\operatorname{Fn}_{<\omega_1}(\omega_1\times\omega_1,2)^L$.

Use the dense suborder of $\mathbb{R}_\lambda$ consisting of the finite functions $f$ assigning pairs in
$[\omega]^{<\omega}\times[\omega_1]^{<\omega}$ to stages in its support.
An omitted stage denotes the greatest condition of its factor.
At a successor one strengthens the preceding condition to decide
the name for the next pair. At a limit a finite support is bounded.
This proves density by induction. The same induction shows that these dense suborders agree in two models $N_0\subseteq N_1$ with the same ordinals
if the factors have the same conditions and order in common generic
extensions by the preceding factors. At a successor, if the two models decided membership or order of a pair differently, this would contradict agreement in a generic extension of the larger model.
At a limit each such decision is made in a bounded initial segment.

The complete embeddings of initial segments and the corresponding translations of names are part of the presentation.
A tail in this presentation has a fixed countable length and uses
Cohen coordinates disjoint from those of the preceding iteration.
The bijections of Cohen index sets used to present products and
quotients belong to $L$.

\begin{lemma}\label{lem:coding-iterations}
The form \eqref{eq:coding-iteration} is preserved by the following operations.
\begin{enumerate}
\item For $\beta\leq\lambda$ and $G_\beta\subseteq \mathbb{P}_\beta$ generic,
\[
 \mathbb{P}_\beta\lessdot \mathbb{P}_\lambda,\qquad
 \mathbb{P}_\lambda/G_\beta\simeq
 \mathbb{C}_{E_{<\lambda}\setminus E_{<\beta}}*\dot {\mathbb{R}}_{[\beta,\lambda)}.
\]
On the dense suborder the projection is
\[
 \pi_\beta(c,f)=
 \bigl(c\restriction(E_{<\beta}\times\omega_1),f\restriction\beta\bigr).
\]
\item A tail defined below $0<d\in\operatorname{RO}(\mathbb{P})$ lifts to a $\mathbb{P}$-name with the same Cohen factors and with
\[
 \dot {\mathbb{A}}_i^{\uparrow}=
 \begin{cases}
  \dot {\mathbb{A}}_i,&\text{on }d,\\
  \mathbf1,&\text{on }\neg d.
 \end{cases}
\]
For $p$ in the dense suborder, $\mathbb{P}_\lambda\restriction p$ again
has form \eqref{eq:coding-iteration}, with initial segments
$\mathbb{P}_\beta\restriction\pi_\beta(p)$, $\beta\leq\lambda$.
\item Tails $\mathbb{Q},\mathbb{R}$ over the same preceding forcing admit product
presentations in either order, on disjoint Cohen index sets.
With $\mathbb{R}$ first, $\mathbb{R}\times \mathbb{Q}_i$ is an initial segment, and
\[
 \mathbb{R}\times \mathbb{Q}_i\lessdot \mathbb{R}\times \mathbb{Q}.
\]
If $\dot {\mathbb{T}}$ is a $\mathbb{Q}_i$-name for a continuation of fixed countable
length, lifted by (2) if necessary, its translation $\dot {\mathbb{T}}_{\mathbb{R}}$ satisfies
\[
 (\mathbb{R}\times \mathbb{Q}_i)*\dot {\mathbb{T}}_{\mathbb{R}}
 \simeq \mathbb{R}\times(\mathbb{Q}_i*\dot {\mathbb{T}})
\]
over $\mathbb{R}\times \mathbb{Q}_i$, preserving all previously specified embeddings
and translations of names.
\end{enumerate}
\end{lemma}

\begin{proof}
(1) Use
\[
 \mathbb{C}_{E_{<\lambda}}\cong
 \mathbb{C}_{E_{<\beta}}\times \mathbb{C}_{E_{<\lambda}\setminus E_{<\beta}}.
\]
If $p=(c,f)$ and $q=(c',f')\leq\pi_\beta(p)$, then
\[
 \bigl(c'\cup c\restriction
       ((E_{<\lambda}\setminus E_{<\beta})\times\omega_1),\,
       f'\cup f\restriction[\beta,\lambda)\bigr)\leq p,q.
\]
After the earlier coordinates are strengthened, the values of $f$ at later stages are still conditions. This proves completeness and the quotient formula,
also at limits.

(2) The lift is obtained by mixing names. For the restriction below $(c,f)$,
choose in $L$, for each $\eta\in\operatorname{supp}(c)$,
\[
 \begin{gathered}
 \sup(\operatorname{dom}c_\eta)+\omega<\nu_\eta<\omega_1,\\
 \tau_\eta:\omega_1\longrightarrow
       \omega_1\setminus\operatorname{dom}c_\eta
       \quad\text{a bijection},\qquad
 \tau_\eta\restriction[\nu_\eta,\omega_1)=\mathrm{id}.
 \end{gathered}
\]
Set $\tau_\eta=\mathrm{id}$ elsewhere. Then
\[
 \psi_c:\mathbb{C}_E\cong \mathbb{C}_E\restriction c,\qquad
 \psi_c(q)_\eta=c_\eta\cup(q_\eta\circ\tau_\eta^{-1}).
\]
Compose the restrictions of $\psi_c$ to the blocks $E_i$ with $\varphi_i$ and replace $\mathbb{A}_i$ by
$\mathbb{A}_i\restriction f(i)$ for $i\in\operatorname{dom}(f)$.
Below $(c,f)$ this preserves $g_i,X_i,Y(X_i)$. By (1) the maps
commute with all initial-segment projections, and a composition of such maps again changes each coordinate only on a bounded set.

(3) Reindex the Cohen coordinates disjointly in $L$ and translate
each tail's names from its own initial segments. In either order
the almost disjoint factors have identical conditions and orders,
so induction on the dense suborders gives the product
isomorphisms. For a continuation of $\mathbb{R}\times \mathbb{Q}_i$, retain both
index maps, extend the $\mathbb{Q}_i$-map, and translate the Boolean value
in (2) by the fixed embedding of $\mathbb{Q}_i$. This gives the displayed
isomorphism over $\mathbb{R}\times \mathbb{Q}_i$. The same induction applies to
translations along an isomorphism and its inverse.
\end{proof}

\begin{lemma}\label{lem:localization-of-reals}
For $\lambda\leq\omega_1$,
\[
 \operatorname{dens}(\mathbb{P}_\lambda)\leq\aleph_1,
 \qquad \mathbb{P}_\lambda\text{ preserves cardinals and }\mathsf{CH}.
\]
If $G_{\omega_1}\subseteq \mathbb{P}_{\omega_1}$ is generic and
$G_\beta=G_{\omega_1}\cap \mathbb{P}_\beta$, then
\[
 (2^\omega)^{W[G_{\omega_1}]}=
 \bigcup_{\beta<\omega_1}(2^\omega)^{W[G_\beta]}.
\]
\end{lemma}

\begin{proof}
Put $\mathbb{C}=\mathbb{C}_{E_{<\lambda}}$ and let $K\subseteq \mathbb{C}$ be generic over $W$.
Since $L\models |\mathbb{C}|\leq\aleph_1$ and $\mathbb{C}$ is countably closed in $L$,
Easton's lemma applied to $\operatorname{Br}$ gives
\[
 (\operatorname{Ord}^{\omega})^{W[K]}
   =(\operatorname{Ord}^{\omega})^W.
\]
Thus $\mathbb{C}$ preserves cardinals and $\mathsf{CH}$ over $W$. In $W[K]$,
$\mathbb{R}_\lambda$ is c.c.c.\ with a dense set of finite functions of size
$\leq\aleph_1$. Hence $\mathbb{C}*\dot {\mathbb{R}}_\lambda$ has the asserted properties.

For $z\in(2^\omega)^{W[G_{\omega_1}]}$, take a nice
$\mathbb{R}_{\omega_1}$-name $\dot z$ in the dense suborder. Its ordinal
code belongs to $W$, and
\[
 \operatorname{supp}(\dot z)=
 \bigcup\{\operatorname{dom}(p):(\check n,p)\in\dot z\}
 \subseteq\beta
 \quad\text{for some }\beta<\omega_1.
\]
The dense suborders of $\mathbb{R}_\beta$ in
$W[K\restriction E_{<\beta}]$ and $W[K]$ agree. Thus $\dot z$ is
already a name over $\mathbb{P}_\beta$, and $z\in W[G_\beta]$.
\end{proof}

\begin{lemma}\label{lem:coded-reals}
In every extension by an iteration of the form \eqref{eq:coding-iteration},
\[
 \{u\in2^\omega:\Phi(u)\}
 =\{u_i:i<\lambda,\ i\text{ is a coding stage}\}.
\]
\end{lemma}

\begin{proof}
The inclusion from right to left follows from
Lemma~\ref{lem:L-decoder}, also for factors restricted below conditions.
Conversely, let $r$ witness $\Phi(u)$. Work first in $W[K]$, with
$K\subseteq \mathbb{C}_{E_{<\lambda}}$ generic, and let $G\subseteq \mathbb{R}_\lambda$
be the remaining generic.

Choose names for $r,u$ and the stage parameters using countable
antichains in the dense suborder. Let $D_i\subseteq i$ collect
the supports of the parameter names at stage $i$, including $a_i$.
Put
\[
 \begin{aligned}
 J_0&=\operatorname{supp}(\dot r)\cup\operatorname{supp}(\dot u),\\
 J_{n+1}&=J_n\cup\bigcup_{i\in J_n}D_i,\qquad
 J=\bigcup_{n<\omega}J_n\in[\lambda]^{\leq\omega}.
 \end{aligned}
\]
The names for $r$, $u$ and for the parameters at stages in $J$ have a common countable ordinal code $z$.
By Easton's lemma and the c.c.c.\ of $\operatorname{Br}$,
\[
 z\in W,\qquad z\in L[H_B]
 \quad\text{for some }B\in[\omega_1]^{\leq\omega}\cap L.
\]
Since $D_i\subseteq J\cap i$ for $i\in J$, induction gives
\[
 \mathbb{R}_{J\cap i}\lessdot \mathbb{R}_i\quad(i\leq\lambda),\qquad
 \dot r^{G_J}=r,\quad\dot u^{G_J}=u,
 \quad G_J=G\cap \mathbb{R}_J.
\]
At successors the factor at a stage in $J$ is given by an $\mathbb{R}_{J\cap i}$-name;
at limits use finite supports and the coherent projections.

Fix $\xi\notin B\cup\bigcup_{i\in J}I_i$ and put
\[
 M_\xi=L[K][H\restriction(\omega_1\setminus\{\xi\})].
\]
Then $z\in M_\xi$. In $M_\xi$ form $\widetilde {\mathbb{R}}_J^\xi$ by following
$\mathbb{R}_J$ until the first stage whose names fail to determine unique
data or whose $I_i$ contains $\xi$, and taking trivial factors
from that stage onward. The test $\xi\in I_i$ precedes the formation
of $X_i$, so the factors before that stage have the same dense suborders over $M_\xi$ and $W[K]$. In $W[K]$, induction gives
\[
 \begin{gathered}
 d_\xi=\bigwedge_{i\in J}
   \|\xi\notin\dot I_i\|_{\operatorname{RO}(\mathbb{R}_J)},\\
 \widetilde {\mathbb{R}}_J^\xi\lessdot \mathbb{R}_J,\qquad
 \widetilde {\mathbb{R}}_J^\xi\restriction d_\xi
   \simeq \mathbb{R}_J\restriction d_\xi.
 \end{gathered}
\]
At limits take meets. Our choice of $\xi$ makes $G_J$ meet
$\mathbb{R}_J\restriction d_\xi$, inducing an $M_\xi$-generic
$\widetilde G\subseteq\widetilde {\mathbb{R}}_J^\xi$ with
\[
 r,u\in M_\xi[\widetilde G].
\]
The Cohen product preserves the c.c.c.\ of ground-model branch
products, whence
\[
 L[K]\models
 \operatorname{Br}_{\omega_1\setminus\{\xi\}}\times\mathbb{S}_\xi
   \text{ is c.c.c.}
\]
All nontrivial factors of $\widetilde {\mathbb{R}}_J^\xi$ are
$\sigma$-centered. Hence
\[
 M_\xi[\widetilde G]\models ``S_\xi\text{ is Suslin}."
\]

Let $\mathcal C=\{i\in J:i\text{ is a coding stage}\}$.
Distinct Cohen coordinates differ on an unbounded set, and the coordinate
isomorphisms change only bounded initial segments. Thus their common
initial segments have bounded length, and
\[
 \sup\bigl(h(g_i)\cap h(g_j)\bigr)<\omega_1
 \qquad(i\ne j\in\mathcal C).
\]
Suppose no stage codes $u$. Choose $\delta<\omega_1$ with
\[
 B\subseteq\delta,\qquad
 h(g_i)\cap h(g_j)\subseteq\delta
 \quad(i\ne j\in\mathcal C).
\]
By Lemma~\ref{lem:L-decoder}, choose $\gamma>\delta$ such that
\[
 \forall k<\omega\quad
 L[r]\models
 ``S_{\omega\gamma+2k+u(k)}\text{ has a cofinal branch}."
\]
There is at most one $i\in\mathcal C$ with $\gamma\in h(g_i)$.
If such $i$ exists, choose $k$ with $u_i(k)\ne u(k)$; otherwise
put $k=0$. Then
\[
 \xi=\omega\gamma+2k+u(k)\notin B\cup\bigcup_{i\in J}I_i.
\]
But $L[r]\subseteq M_\xi[\widetilde G]$, and $S_\xi$ is Suslin in that model, a contradiction.
\end{proof}

Up to the specified presentations, the bounded iterations, their
restrictions below conditions, and the names for their tails form
a set. Choose dense suborders carried by subsets of $\omega_1$.
Their orders, the specified embeddings, and the names for stage
parameters can then be encoded by subsets of fixed sets of conditions
and names. For a tail whose countable length is given by a name,
first work below a condition deciding the length and apply
Lemma~\ref{lem:coding-iterations}(2). Thus the forcing predicates
used below may be indexed by tail names of fixed countable length.

\subsection{Motivation for the allowability hierarchy}
\label{subsection-motivation-allowability}

We call the coding iterations of countable length from Subsection~\ref{subsec:coding-iterations}, possibly restricted below a condition, and their quotients over extensions by initial segments \emph{allowable}.  Given the notion of allowable forcing, we now form successive approximations to the class of forcings which we eventually want to use.  We first describe the idea leading to these approximations.

Fix a recursive enumeration
\[
   \langle \varphi_m(v_0,v_1)\mid m\in\omega\rangle
\]
of the $\Pi^1_3$-formulas with two free real variables, and put
\[
   A_m:=\{(x,y)\in(2^\omega)^2:\varphi_m(x,y)\}.
\]
We write $f_m$ for the uniformizing function which we intend to define for $A_m$, and write $f(m,x)$ for $f_m(x)$.  The goal is to choose, for every $m\in\omega$ and every real $x$ for which
\[
   A_{m,x}:=\{y\in2^\omega:(x,y)\in A_m\}
\]
is nonempty, a real $f(m,x)$ such that
\[
   (x,f(m,x))\in A_m,
\]
and to code every triple $(x,y,m)$ such that
\[
   y\neq f(m,x)\quad\text{and}\quad(x,y)\in A_m
\]
into $\vec S$.  The final forcing will have length $\omega_1$, and each of its countable initial segments will be allowable.  Since the statement
\[
   ``(x,y,m)\text{ is coded into }\vec S"
\]
is $\Sigma^1_3$, the requirement
\[
   (x,y)\in A_m
   \quad\text{and}\quad
   (x,y,m)\text{ is not coded into }\vec S
\]
is $\Pi^1_3$.  The bookkeeping of the final iteration will ensure that, whenever $A_{m,x}\neq\varnothing$, exactly one such $y$ remains.

The idea behind the definition can be seen in the following toy example.  Work in $W$, let $x\in W$, and suppose that
\[
   A_{m,x}=\{y_0,y_1\}.
\]
Assume in addition that no allowable forcing over $W$ adds a new member to $A_{m,x}$.  We consider only the problem of defining $f(m,x)$ and omit, for the moment, the interaction with the other sections.

Call $y\in A_{m,x}$ \emph{eliminable} if there is an allowable forcing $\forceP$ such that
\[
   \forceP\Vdash (x,y)\notin A_m.
\]
Call $y$ \emph{persistent} if it is not eliminable.  Thus, if $y_0$ is persistent, we may reserve $y_0$ as the value of $f(m,x)$ and code $(x,y_1,m)$ into $\vec S$.  If $y_0$ is eliminable, witnessed by an allowable forcing $\forceP$, we consider $y_1$.  If $y_1$ is persistent, we reserve $y_1$ and code $(x,y_0,m)$.

It remains to consider the case in which both values are eliminable.  Let $\forceP$ and $\forceQ$ witness, respectively, that
\[
   \forceP\Vdash (x,y_0)\notin A_m
   \qquad\text{and}\qquad
   \forceQ\Vdash (x,y_1)\notin A_m.
\]
By Lemma~\ref{lem:coding-iterations}, the product
\[
   \forceP\times\forceQ
\]
is isomorphic to an allowable forcing.  Using this presentation of the product, upward absoluteness gives
\[
   \forceP\times\forceQ\Vdash
   (x,y_0)\notin A_m\wedge(x,y_1)\notin A_m.
\]
By the assumption that no new member is added to the section, the product forces
\[
   A_{m,x}=\varnothing.
\]
In this case there is no value to uniformize at $x$.

This argument isolates the two operations which are needed later:
\begin{equation}
\label{two-operations-uniformization}
\begin{split}
   &\text{reserve a value which is persistent under the forcings under consideration;}\\
   &\text{combine eliminating forcings by taking products.}
\end{split}
\end{equation}
The difficulty is that an eliminating forcing may add codes which interfere with other sections.  Thus it is not enough to use arbitrary allowable forcings as witnesses.  We need a class $\mathcal P$ of allowable forcings such that the following rule produces forcings in $\mathcal P$ again.  For every task $(x,m)$ considered by the bookkeeping, either
\[
   \exists y\in A_{m,x}\ \forall\forceQ\in\mathcal P
   \bigl(\forceQ\text{ extending the current forcing}\rightarrow
   \forceQ\Vdash(x,y)\in A_m\bigr),
\]
in which case one such $y$ is reserved and the other values are coded, or every value under consideration can be eliminated by a forcing in $\mathcal P$.

This turns the construction into a fixed point problem.  Let
\[
   \mathcal A^W_0
\]
be the class of allowable iterations.  The first refinement searches for values which are persistent under $\mathcal A^W_0$-extensions; the resulting iterations form $\mathcal A^W_1$.  The next refinement searches also for values which are persistent under $\mathcal A^W_1$-extensions, and so on.  We obtain a decreasing hierarchy
\[
   \mathcal A^W_0\supseteq\mathcal A^W_1\supseteq\cdots
   \supseteq\mathcal A^W_\alpha\supseteq\cdots.
\]
In the first refinement we do not code a minimizer if all values of least rank are already present.

The reason for continuing the construction is the usual moving-target problem.  In the definition of $1$-allowability we test persistence under $0$-allowable extensions, but the iteration which we produce is itself $1$-allowable.  We should therefore also test persistence under $1$-allowable extensions.  Repeating this argument leads to the transfinite hierarchy defined in the next subsection.

\subsection{Persistence and the choice of a uniformizing value}
\label{subsec:persistence}\label{subsection-alpha-allowability}

Use the coding iterations of Subsection~\ref{subsec:coding-iterations},
with the specified complete embeddings of their initial segments.
We work in extensions \(M\) by such an initial segment, sometimes
restricted below a condition. An iteration over \(M\) is a quotient
with the presentation given in that subsection. The phrase
\emph{code \((m,x,y)\)} means to code the real \(\langle x,y,m\rangle\).
Put
\[
   A(m,x,y)\quad\Longleftrightarrow\quad (x,y)\in A_m .
\]
We use two absoluteness facts for \(A\). A witness to \(\neg A(b)\)
in the extension by one factor of a product remains a witness in
the product extension. If \(\neg A(b)\) holds in the final extension,
its witness gives \(\neg A(b)\) in the extension by every sufficiently long proper initial segment of the iteration. Here \(b=(m,x,y)\) belongs to the
model in which membership is tested. In the construction over \(L\),
write \(\neg A(b)\) as \(\exists z\,\psi(b,z)\) with
\(\psi\in\Pi^1_2\) and apply Shoenfield absoluteness to \(\psi\).
The second assertion also uses the fact that every real belongs to
an extension by a bounded initial segment.

\paragraph{Choice functions.}
For each extension \(M\) under consideration we include a name for
a choice function \(J_M\) on its nonempty sets of reals, so
\(J_M(B)\in B\). Along each specified inclusion \(M\subseteq N\),
require \(J_N(B)=J_M(B)\) for \(B\in\operatorname{dom}(J_M)\).
These are names over the given forcings.

\begin{lemma}
\label{lem:choice-functions}
Let \(\mathbb{Q},\mathbb{R}\) be countable iterations over \(M\). For their product
extensions put
\[
 M_{ij}=M[G_{\mathbb{Q}_i},G_{\mathbb{R}_j}],\qquad
 \mathcal D_{ij}=
   \bigl(\mathcal P(\mathbb R)\setminus\{\varnothing\}\bigr)^{M_{ij}}.
\]
Coherent choice functions on the extensions \(M_{i0}\) and \(M_{0j}\), agreeing on \(M=M_{00}\), extend by names to functions
\[
 \begin{gathered}
 J_{ij}:\mathcal D_{ij}\longrightarrow\mathbb R^{M_{ij}},
 \qquad J_{ij}(B)\in B,\\
 i\leq i',\ j\leq j'
 \quad\Longrightarrow\quad
 J_{i'j'}\restriction\mathcal D_{ij}=J_{ij}.
 \end{gathered}
\]
If coherent functions are prescribed for \(i\leq i_0,\ j\leq j_0\), the extension can be chosen to agree with them. The assertions remain
valid under restrictions below conditions, mixing of names,
quotients by \(\mathbb{Q}_{i_0}\times \mathbb{R}_{j_0}\), and appending an iteration
to either factor of a specified product.
\end{lemma}

\begin{proof}
Order pairs coordinatewise, with meet \(a\wedge b\) and join
\(a\vee b\). We first show
\[
 \mathcal D_a\cap\mathcal D_b=\mathcal D_{a\wedge b}.
\]
For incomparable \(a,b\), factor over \(N=M_{a\wedge b}\).
Fix a bijection \(e:\theta\to\mathbb R^N\) in \(N\).
The intersection lemma for the mutually generic quotients, applied first
to each element of \(B\) and then to \(e^{-1}[B]\), gives
\[
 B\in\mathcal D_a\cap\mathcal D_b
 \ \Longrightarrow\ B\subseteq\mathbb R^N
 \ \Longrightarrow\ e^{-1}[B]\in N
 \ \Longrightarrow\ B\in\mathcal D_{a\wedge b}.
\]
The comparable case is immediate. Consequently
\[
 \ell(B)=\bigl(\min\{i:\exists j\ B\in\mathcal D_{ij}\},\,
                  \min\{j:\exists i\ B\in\mathcal D_{ij}\}\bigr)
\]
satisfies
\[
 B\in\mathcal D_{\ell(B)},\qquad
 B\in\mathcal D_a\ \Longleftrightarrow\ \ell(B)\leq a.
\]
Indeed, take the meet of pairs witnessing the two coordinate minima.

Let \(J_a^0\) denote the prescribed functions and let
\(\triangleleft_a\) be a named well-order of \(\mathbb R^{M_a}\).
For \(B\in\mathcal D_a\), set
\[
 J_a(B)=
 \begin{cases}
  J^0_{\ell(B)}(B),&\text{if }\ell(B)\text{ is prescribed},\\
  \min_{\triangleleft_{\ell(B)}}B,&\text{otherwise}.
 \end{cases}
\]
The set of pairs with prescribed functions is downward closed, so all prescribed values are kept and \(a\leq b\Rightarrow
J_b\restriction\mathcal D_a=J_a\).
These are names: for \(a\leq b\), membership in the smaller
extension is expressed by
\[
 \|\dot B\in M_a\|
   =\bigvee_{\dot C}\|\dot B=\iota_{ab}(\dot C)\|,
\]
where \(\iota_{ab}\) is the specified translation and \(\dot C\)
ranges over a fixed set of names for sets of reals over
\(\mathbb{Q}_{a_0}\times \mathbb{R}_{a_1}\), where \(a=(a_0,a_1)\).
For \(B\in\mathcal D_b\), the minima defining \(\ell(B)\)
can be computed using only pairs \(a\leq b\).

If \(B\in\mathcal D_c\) is old and occurs at a new pair \(a\), then
\(B\in\mathcal D_{a\wedge c}\), with \(a\wedge c\) old. Hence
\(\ell_{\mathrm{new}}(B)=\ell_{\mathrm{old}}(B)\).
For \(s=(i_0,j_0)\), restricting to pairs \(a\geq s\) instead gives
\[
 \ell_{\geq s}(B)=\ell(B)\vee s.
\]
If \(\ell(B)\vee s\ne\ell(B)\), then \(\ell(B)\vee s\) has first coordinate \(i_0\) or second coordinate \(j_0\), and the choice is that of the original product; otherwise \(\ell(B)\vee s=\ell(B)\), so the original choice applies.
This gives compatibility with quotients by
\(\mathbb{Q}_{i_0}\times \mathbb{R}_{j_0}\) and with both product orders.

For restrictions use the specified embeddings and
\[
 (\mathbb{Q}\times \mathbb{R})\restriction(q,r)
   =(\mathbb{Q}\restriction q)\times(\mathbb{R}\restriction r).
\]
Below \(d\) retain the prescribed names; on \(\neg d\) use a
coherent extension and mix the names. Finally, for a continuation
\(\dot {\mathbb{T}}\) over \(\mathbb{Q}_{i_0}\) given below \(d\),
Lemma~\ref{lem:coding-iterations}(2) supplies \(\dot {\mathbb{T}}^{\uparrow}\)
with trivial almost disjoint factors off \(d\). Apply the construction
to \(\mathbb{Q}_{i_0}*\dot {\mathbb{T}}^{\uparrow}\) and the other factor.
Since \(\ell(B)\) does not change, the earlier choices are kept, and Lemma~\ref{lem:coding-iterations}(3) preserves the specified embeddings.
\end{proof}

We work with the set of all countable iterations under consideration,
together with their specified embeddings and the names for their
choice functions. The same forcing may occur with different data.
For each presentation \(\mathbb H\), fix a dense suborder
\(D_{\mathbb H}\) of size at most \(\aleph_1\). Names for reals may
be taken in the form
\[
 \dot r=\{(\check k,p):k<\omega,\ p\in A_k\},
 \qquad A_k\subseteq D_{\mathbb H}\text{ an antichain}.
\]
These names form a fixed set. Names for sets of reals are represented
by families of antichains indexed by those real names; the same
construction, applied to graphs, supplies names for choice functions.
Thus the presentations together with the names and embeddings still
form a set, possibly of size greater than \(2^{\aleph_1}\).
This set includes the restrictions and products in
Lemma~\ref{lem:choice-functions}. A name for a quotient iteration
defined below a condition can be used there or extended by trivial
almost disjoint factors on the complement, with the same countable
length.

\paragraph{The allowability hierarchy.}
We define successive hierarchies of allowable iterations. Each
hierarchy begins at a specified initial segment of the final
iteration. The \(\xi\)-th hierarchy keeps the restrictions on
coding imposed so far and the fixed classes from hierarchies
\(\eta<\xi\). A restriction forbids coding a specified triple
\((m,x,y)\). At a later extension by a bounded initial segment,
evaluate the earlier definitions via the complete embeddings from
the initial segments at which those definitions began. In a
product with a further iteration these embeddings are obtained
through the factor belonging to the final iteration. 

For a generic extension \(M\) considered in the current hierarchy, let the current class of base forcings
\(\mathcal K^\xi_0(M)\) consist of the quotient iterations over \(M\)
that obey all the restrictions on coding and belong to every
earlier fixed class, evaluated over \(M\).  Assume now that we have defined \(\mathcal K^\xi_\beta\) in all these
extensions, put
\[
 P^\xi_\beta(M)=
 \{b\in(\omega\times\mathbb R^2)^M:
           \forall \mathbb{Q}\in\mathcal K^\xi_\beta(M)\quad M \models \mathbb{Q}\Vdash A(b)\}.
\]
The rank of \(b\) up to \((\xi,\beta)\) is the first pair, in
lexicographic order, at which \(b\) belongs to a persistence set.
Use all levels of the earlier hierarchies up to their fixed
levels and the levels at most \(\beta\) of the current hierarchy.
Write \(\rho_{\leq(\xi,\beta)}(b)\) for this partial rank.
It is undefined precisely when \(b\) has no rank at any of these
levels, including those in the earlier hierarchies.

For \(m,x\in M\), let \(\Lambda^\xi_\beta(M,m,x)\) be the set of
pairs \(\lambda\) for which some \(\mathbb{Q}\in\mathcal K^\xi_\beta(M)\),
some \(q\in \mathbb{Q}\), and some real name \(\dot y\) satisfy
\[
 q\Vdash_{\mathbb{Q}}^M \rho_{\leq(\xi,\beta)}(m,x,\dot y)=\lambda.
\]
When this set is nonempty, put
\(\lambda=\min\Lambda^\xi_\beta(M,m,x)\). Call this minimum
\emph{stable} if no condition witnessing it forces
\(\dot y\notin\mathbb R^M\). In this case let \(B\) be the
nonempty set of \(y\in\mathbb R^M\) for which some such condition
forces \(\rho_{\leq(\xi,\beta)}(m,x,\check y)=\lambda\), and put
\(j^\xi_\beta(M,m,x)=J_M(B)\). Thus \(j^\xi_\beta\) is defined
precisely in the stable case.

An iteration over \(M\) belongs to \(\mathcal K^\xi_{\beta+1}(M)\)
if it belongs to \(\mathcal K^\xi_\beta(M)\) and, whenever a stage
codes \((m,x,y)\), the extension \(N\) by the preceding initial
segment satisfies
\begin{equation}\label{eq:requirement}
       j^\xi_\beta(N,m,x)\downarrow
        \quad\Longrightarrow\quad y\ne j^\xi_\beta(N,m,x).
\end{equation}
Take intersections at limit levels. At a limit of hierarchies
retain the union of the restrictions on coding and use every
earlier fixed class in the definition at level zero. We will always be able to produce non-empty classes which form a fixed point in the following way:

\begin{lemma}
\label{lem:generic-recursion}
Let \(M_0\) be the model at the start of hierarchy \(\xi\).
The preceding recursion uniformly defines \(\mathcal K^\xi_\beta(M)\)
in every generic extension \(M=M_0[G]\) under consideration.
There is an ordinal \(\theta_\xi\), independent of \(M\),
such that
\[
 \beta\geq\theta_\xi
 \quad\Longrightarrow\quad
 \mathcal K^\xi_\beta(M)=\mathcal K^\xi_{\theta_\xi}(M)
 \qquad\text{for every such }M.
\]
\end{lemma}

\begin{proof}
All hypothetical generic filters below are interpreted by the forcing
theorem. Induct on \(\beta\), simultaneously over the chosen forcing
presentations. Let \(G\subseteq\mathbb H\) be \(M_0\)-generic
and put \(M=M_0[G]\). At level zero, retain the restrictions on
coding and the earlier fixed classes, interpreted through the
specified embeddings.

Suppose level \(\beta\) is defined, and abbreviate
\(\rho^N_{\leq(\xi,\beta)}\) by \(\rho^N_\beta\).
For \(b\in(\omega\times\mathbb R^2)^M\), \(m\in\omega\), and
\(x\in\mathbb R^M\), the definitions say
\[
\begin{aligned}
 b\in P^\xi_\beta(M)
 &\ \Longleftrightarrow\
 \forall\mathbb Q,H\quad
 M[H]\models A(b),\\
 \lambda\in\Lambda^\xi_\beta(M,m,x)
 &\ \Longleftrightarrow\
 \exists\mathbb Q,H,y\quad
 \rho^{M[H]}_\beta(m,x,y)=\lambda,
\end{aligned}
\]
where \(\mathbb Q\in\mathcal K^\xi_\beta(M)\),
\(H\subseteq\mathbb Q\) ranges over hypothetical \(M\)-generic
filters, and \(y\in\mathbb R^{M[H]}\). At the least attainable \(\lambda\),
stability means that no witness \(y\) in the second clause lies
outside \(M\). In the stable case, let \(B\subseteq\mathbb R^M\)
be the nonempty set of witnesses from that clause and put
\(j^\xi_\beta(M,m,x)=J_M(B)\).

For a candidate iteration \(\mathbb T\) over \(M\), let
\(H\subseteq\mathbb T\) be hypothetical \(M\)-generic. At a coding
stage \(i\), write \(M_i=M[H_i]\) for the preceding intermediate
extension and \((m_i,x_i,y_i)\) for the triple coded there.
Then \(\mathbb T\in\mathcal K^\xi_{\beta+1}(M)\) iff
\(\mathbb T\in\mathcal K^\xi_\beta(M)\) and, for every such \(H\)
and every coding stage \(i\),
\[
 j^\xi_\beta(M_i,m_i,x_i)\downarrow
 \quad\Longrightarrow\quad y_i\ne j^\xi_\beta(M_i,m_i,x_i).
\]
At a limit \(\delta\),
\(\mathcal K^\xi_\delta(M)=\bigcap_{\beta<\delta}\mathcal K^\xi_\beta(M)\).
The induction hypothesis and the forcing theorem express every clause
uniformly by conditions and names in the preceding model.

When restricting or factoring the forcing, use the usual translation
of names for all parameters. At a limit of hierarchies, intersect the
earlier fixed classes.

Finally, in \(M_0\) fix the set \(X_\xi\) of triples
\((\mathbb H,\dot{\mathbb T},p)\), where \(\mathbb H,\dot{\mathbb T}\)
are among the chosen presentations and \(p\in\mathbb H\), and put
\[
 D_\beta=\{(\mathbb H,\dot{\mathbb T},p)\in X_\xi:
            p\Vdash_{\mathbb H}\dot{\mathbb T}\in\mathcal K^\xi_\beta\}.
\]
Here the forcing statement uses the interpretation just defined.
For \(\beta\leq\gamma\), \(D_\gamma\subseteq D_\beta\).
The least stages at which members of \(D_0\) disappear form a set;
choose \(\theta_\xi\) above all these stages. Then
\(D_\beta=D_{\theta_\xi}\) for every \(\beta\geq\theta_\xi\).
For every \(M_0\)-generic \(G\subseteq\mathbb H\), the truth lemma gives
\[
 \dot{\mathbb T}^{G}\in\mathcal K^\xi_\beta(M_0[G])
 \quad\Longleftrightarrow\quad
 \exists p\in G\ (\mathbb H,\dot{\mathbb T},p)\in D_\beta.
\]
Since the chosen names cover all relevant tails, the same
\(\theta_\xi\) works in every such generic extension. Put
\(K_\xi(M)=\mathcal K^\xi_{\theta_\xi}(M)\).
\end{proof}

\begin{lemma}
\label{lem:products-ranks}
If \(\mathbb{Q},\mathbb{R}\in\mathcal K^\xi_\alpha(M)\), then
\[
 (\mathbb{Q}\times \mathbb{R})/G_{\mathbb{Q}}\in\mathcal K^\xi_\alpha(M[G_{\mathbb{Q}}]),
 \qquad
 (\mathbb{Q}\times \mathbb{R})/G_{\mathbb{R}}\in\mathcal K^\xi_\alpha(M[G_{\mathbb{R}}]).
\]
For \(N=M[G_{\mathbb{R}}]\), write \(\rho^M,\rho^N\) for the partial ranks up to \((\xi,\alpha)\). Then
\[
 \rho^N\restriction(\omega\times\mathbb R^2)^M=\rho^M,
 \qquad
 \Lambda^\xi_\alpha(N,m,x)=\Lambda^\xi_\alpha(M,m,x)
 \quad(m,x\in M).
\]
If the common minimum exists, it is stable over \(M\) iff it is
stable over \(N\). In the stable case the sets of minimizers agree and
\(j^\xi_\alpha(N,m,x)=j^\xi_\alpha(M,m,x)\).
\end{lemma}

\begin{proof}
The proof proceeds  by induction lexicographically on \((\xi,\alpha)\), proving the product
assertions first. Evaluating the defining requirements in the intermediate
generic extensions, as in Lemma~\ref{lem:generic-recursion}, gives closure under
\[
 \mathbb{T}\mapsto \mathbb{T}\restriction t,\qquad
 \mathbb{T}\mapsto \mathbb{T}_i,\qquad
 \mathbb{T}\mapsto \mathbb{T}/G_{\mathbb{T}_i}.
\]
A name given below \(d\) is lifted using trivial almost disjoint
factors on \(\neg d\). Evaluating each requirement immediately
before its coding stage also gives
\[
 \begin{gathered}
 \mathbb{T}\in\mathcal K^\xi_\alpha(S),\qquad
 \mathbb{T}\Vdash\dot {\mathbb{U}}\in\mathcal K^\xi_\alpha(S[G_{\mathbb{T}}])\\
 \Longrightarrow\quad \mathbb{T}*\dot {\mathbb{U}}\in\mathcal K^\xi_\alpha(S).
 \end{gathered}
\]

In \(\mathbb{Q}\times \mathbb{R}\), each coding stage codes the same real as in its factor, so the coding restrictions hold. Membership in earlier fixed
classes follows from the induction on \(\xi\), through the
specified embeddings. For \(\beta<\alpha\), let \(S\) be the
extension preceding a stage of one factor and \(S'=S[G_{\mathbb{U}}]\)
the extension obtained by adding the other factor. By induction,
\[
 \mathbb{U}\in\mathcal K^\xi_\beta(S),\qquad
 j^\xi_\beta(S',m,x)=j^\xi_\beta(S,m,x),
\]
where equality includes definedness. Thus
\eqref{eq:requirement} transfers from \(S\) to \(S'\), proving
the product assertions at \(\alpha\).

Fix \(N=M[G_{\mathbb{R}}]\) and \(b\in(\omega\times\mathbb R^2)^M\). We have
\[
 b\notin P^\xi_\alpha(N)
 \quad\Longleftrightarrow\quad
 \exists \mathbb{T}\in\mathcal K^\xi_\alpha(N)\ \exists t\in \mathbb{T}\
 t\Vdash_{\mathbb{T}}^N\neg A(b).
\]
Given such \(\mathbb{T},t\), choose their names and decide the length below
\(r\in G_{\mathbb{R}}\). Lift the iteration as above. Then
\[
 \mathbb{R}*\dot {\mathbb{T}}^\uparrow\in\mathcal K^\xi_\alpha(M),\qquad
 (r,\dot t)\Vdash\neg A(b),
\]
so \(b\notin P^\xi_\alpha(M)\). Conversely, if
\(\mathbb{Q}\in\mathcal K^\xi_\alpha(M)\) and
\(q\Vdash_{\mathbb{Q}}^M\neg A(b)\), the witness persists in \(\mathbb{Q}\times \mathbb{R}\).
Hence, identifying \(q\) with its quotient condition,
\[
 q\Vdash^N_{(\mathbb{Q}\times \mathbb{R})/G_{\mathbb{R}}}\neg A(b).
\]
Consequently
\[
 P^\xi_\alpha(N)\cap M=P^\xi_\alpha(M).
\]
Applying this at the lower levels and earlier hierarchies yields
\(\rho^N\restriction(\omega\times\mathbb R^2)^M=\rho^M\).

Composition \(\mathbb{R}*\dot {\mathbb{T}}^\uparrow\) and the product \(\mathbb{T}\times \mathbb{R}\)
give, respectively, the two inclusions in
\[
 \Lambda^\xi_\alpha(N,m,x)=\Lambda^\xi_\alpha(M,m,x).
\]
For the product direction, the quotient over \(M[G_{\mathbb{T}}]\) belongs to \(\mathcal K^\xi_\alpha(M[G_{\mathbb{T}}])\), so the witnessing rank is preserved. For new witnesses, the intersection lemma for mutually generic extensions gives
\[
 \begin{gathered}
 M[G_{\mathbb{T}}]\cap M[G_{\mathbb{R}}]=M,\\
 y\in M[G_{\mathbb{T}}]\setminus M\quad\Longrightarrow\quad y\notin N.
 \end{gathered}
\]
Conversely, \(y\notin N\Rightarrow y\notin M\) in a composed
extension. Thus the common minimum is stable over \(M\) iff it
is stable over \(N\).

In the stable case, put \(\lambda=\min\Lambda^\xi_\alpha(M,m,x)\)
and let \(B_M,B_N\) be the minimizer sets. Composition gives \(B_N\subseteq B_M\), and the product gives
\(B_M\subseteq B_N\). Rank preservation and the trivial forcing give
\[
 B_S=\{y\in\mathbb R^S:\rho^S(m,x,y)=\lambda\}
 \qquad(S=M,N).
\]
For \(B=B_M=B_N\), coherence therefore yields
\[
 j^\xi_\alpha(N,m,x)=J_N(B)=J_M(B)
                   =j^\xi_\alpha(M,m,x).
\]
\end{proof}

In particular, whenever \eqref{eq:requirement} from an
earlier level or hierarchy forbids coding a triple over \(M\),
the minimum remains stable and the same triple is forbidden
by the corresponding requirement at the current fixed level.
Indeed, a minimizer already in \(M\) witnesses \(\lambda\)
using the trivial forcing in the smaller class. Any smaller
rank attainable in the smaller class would already have been attained at the earlier level, since newly assigned ranks are
later in the lexicographic order. Any new real of rank
\(\lambda\) obtainable in the smaller class would also have
been obtainable in the earlier class. Thus the set of
minimizers and its choice agree. If the current minimum can be attained by a new real, no earlier
requirement of this kind is defined at \((m,x)\) over \(M\).

\begin{lemma}
\label{lem:new-restrictions}
Let \(M\) be an extension by a bounded initial segment of the
final iteration during hierarchy \(\eta\), and let \(\mathbb{R}\in K_\eta(M)\)
be a further bounded segment of the final iteration. Suppose
\[
 (m,x,y)\text{ forbidden during }\mathbb{R}\quad\Longrightarrow\quad y\notin M .
\]
Let \(\xi\) be the current hierarchy in \(M[G_{\mathbb{R}}]\). For every
\(\mathbb{Q}\in K_\eta(M)\), the quotient \((\mathbb{Q}\times \mathbb{R})/G_{\mathbb{R}}\) belongs to all
the intervening fixed classes and to \(K_\xi(M[G_{\mathbb{R}}])\).
For each initial segment \(\mathbb{Q}_i\), put
\[
 S=M[G_{\mathbb{Q}_i}],\qquad N=M[G_{\mathbb{Q}_i},G_{\mathbb{R}}].
\]
Writing \(\rho^S,\rho^N\) for the partial ranks up to \((\eta,\theta_\eta)\) and \((\xi,\theta_\xi)\), respectively,
we have
\[
 \rho^N\restriction(\omega\times\mathbb R^2)^S=\rho^S.
\]
For \((m,x)\in S\), if the minimum over \(S\) in \(K_\eta(S)\)
exists, then
\[
 \min\Lambda^\xi_{\theta_\xi}(N,m,x)
   =\min\Lambda^\eta_{\theta_\eta}(S,m,x).
\]
The two minima have the same stability. In the stable case the
sets of minimizers and the choices agree.
\end{lemma}

\begin{proof}
Fix the product presentation and induct on the later hierarchies
and their levels, simultaneously for all \(\mathbb{Q}\in K_\eta(M)\) and
all \(\mathbb{Q}_i\). In each bounded extension containing the start
of an earlier hierarchy, interpret its definitions using the generic
filters induced through the final-iteration factor. Retain this
interpretation after a new hierarchy has started.
Lemma~\ref{lem:coding-iterations}(3) retains the embeddings under
\(\mathbb{Q}_i\mapsto\mathbb{Q}_i*\dot{\mathbb{T}}\), and
Lemma~\ref{lem:choice-functions} extends the choice functions coherently.

At level zero, the earlier fixed classes follow from
Lemma~\ref{lem:products-ranks} and induction on the hierarchies.
A newly forbidden \((m,x,y)\) cannot be coded by \(\mathbb{Q}\): at its
preceding stage \(\mathbb{Q}_i\) we would have
\[
 y\in M[G_{\mathbb{Q}_i}]\cap M[G_{\mathbb{R}}]=M,\qquad y\notin M.
\]
At a limit of hierarchies intersect the earlier fixed classes.

Suppose the product assertion holds at level \(\beta\) of
hierarchy \(\xi\). Put
\[
 S=M[G_{\mathbb{Q}_i}],\qquad N=M[G_{\mathbb{Q}_i},G_{\mathbb{R}}],\qquad
 \rho_\beta^N=\rho_{\leq(\xi,\beta)}^N.
\]
For \(b\in(\omega\times\mathbb R^2)^S\),
\[
 \rho^S(b)\uparrow
 \quad\Longrightarrow\quad
 \exists \mathbb{T}\in K_\eta(S)\ \exists t\in \mathbb{T}\
 t\Vdash_{\mathbb{T}}^S\neg A(b).
\]
Replace \(\mathbb{T}\) by \(\mathbb{T}\restriction t\), choose its name and length
below a condition in \(G_{\mathbb{Q}_i}\), and lift it using trivial almost
disjoint factors on the complement. Then
\[
 \mathbb{Q}_i*\dot {\mathbb{T}}^\uparrow\in K_\eta(M).
\]
Extend the specified presentation of \(\mathbb{Q}_i\times \mathbb{R}\) to
\((\mathbb{Q}_i*\dot {\mathbb{T}}^\uparrow)\times \mathbb{R}\). By induction its appended quotient
\[
 \mathbb{U}=\bigl((\mathbb{Q}_i*\dot {\mathbb{T}}^\uparrow)\times \mathbb{R}\bigr)/
                    (G_{\mathbb{Q}_i}\times G_{\mathbb{R}})
\]
satisfies
\[
 \mathbb{U}\in\mathcal K^\xi_\beta(N),\qquad
 \mathbb{U}\Vdash_N\neg A(b),\qquad b\notin P^\xi_\beta(N).
\]
Here the witness to \(\neg A(b)\) persists in the product.
Previously assigned ranks agree by
Lemma~\ref{lem:products-ranks}, evaluating the earlier
definitions in this same product. Together with the induction,
\[
 \rho_\beta^N\restriction(\omega\times\mathbb R^2)^S=\rho^S.
\]

Suppose
\(\lambda=\min\Lambda^\eta_{\theta_\eta}(S,m,x)\) exists.
Let \(\mathbb{T}\in K_\eta(S)\) witness \(\lambda\), restricted below a
witnessing condition. Appending \(\mathbb{T}\) as above gives
\[
 \lambda\in\Lambda^\xi_\beta(N,m,x).
\]
In the product over \(S\), the intersection lemma gives
\[
 S[G_{\mathbb{T}}]\cap N=S,\qquad
 y\in S[G_{\mathbb{T}}]\setminus S\quad\Longrightarrow\quad y\notin N.
\]
Thus a new witness remains new. Conversely, put
\(\mathbb{V}=(\mathbb{Q}_i\times \mathbb{R})/G_{\mathbb{Q}_i}\). By
Lemma~\ref{lem:products-ranks}, \(\mathbb{V}\in K_\eta(S)\).
For a witnessing \(\mathbb{U}\in\mathcal K^\xi_\beta(N)\), choose and lift
its name as above. Then
\[
 \mathbb{V}*\dot {\mathbb{U}}^\uparrow\in K_\eta(S).
\]
Since every rank \(\mu<\lambda\) is already a rank up to \((\eta,\theta_\eta)\), this composition gives
\[
 \mu<\lambda,\quad\mu\in\Lambda^\xi_\beta(N,m,x)
 \quad\Longrightarrow\quad
 \mu\in\Lambda^\eta_{\theta_\eta}(S,m,x),
\]
contrary to minimality. Hence
\[
 \min\Lambda^\xi_\beta(N,m,x)=\lambda.
\]
If the old minimum is stable, composition excludes new minimizers
over \(S\); the old minimizers retain rank \(\lambda\).
For the corresponding sets \(B_S,B_N\), therefore,
\[
 B_N\subseteq B_S\subseteq B_N,\qquad
 j^\xi_\beta(N,m,x)=J_N(B_N)=J_S(B_S)
                 =j^\eta_{\theta_\eta}(S,m,x).
\]
Together with the preservation of new witnesses, this proves
equality of stability.

At a stage of \(\mathbb{Q}\) coding \((m,x,y)\), if the old minimum exists,
its stability and choice transfer \eqref{eq:requirement}
from \(S\) to \(N\). Otherwise, the trivial forcing and rank
preservation give
\[
 \Lambda^\eta_{\theta_\eta}(S,m,x)=\varnothing
 \ \Longrightarrow\ \rho^S(m,x,y)\uparrow
 \ \Longrightarrow\ \rho_\beta^N(m,x,y)\uparrow.
\]
Thus \(j^\xi_\beta(N,m,x)\downarrow\) implies
\(y\ne j^\xi_\beta(N,m,x)\) in this case as well.
This proves the successor step. At limit levels take all earlier
requirements, yielding the assertions at the fixed level.
\end{proof}

\begin{lemma}
\label{lem:first-segment-minimizer}
Suppose the current fixed class allows its minimum rank
\(\lambda\) at \((m,x)\) to be attained by a real outside \(M\).
There is a restriction of an allowable iteration that forces
a real \(y\) of rank \(\lambda\) at \((m,x)\) to appear in its
final extension and in no extension by a proper initial segment.
\end{lemma}

\begin{proof}
Choose \(\mathbb{Q}\in K_\xi(M)\), \(q\in \mathbb{Q}\), and \(\dot y\) with
\[
 q\Vdash_{\mathbb{Q}}
 \dot y\notin\mathbb R^M\ \land\
 \rho_{\leq(\xi,\theta_\xi)}(m,x,\dot y)=\lambda.
\]
Strengthen \(q\) to decide \(\delta>0\) and a \(\mathbb{Q}_\delta\)-name
\(\dot z\) such that
\[
 \begin{aligned}
 q&\Vdash_{\mathbb{Q}}
 \delta=\min\{\gamma\leq\operatorname{lh}(\mathbb{Q}):
                   \dot y\in\dot{\mathbb R}^{\mathbb{Q}_\gamma}\},\\
 q&\Vdash_{\mathbb{Q}}\dot y=\dot z,
 \end{aligned}
\]
identifying names with their lifts.
Since \(\mathbb{Q}/G_{\mathbb{Q}_\delta}\) is allowable,
Lemma~\ref{lem:products-ranks} implies that \(q\) forces
the \(\mathbb{Q}_\delta\)-extension to satisfy
\[
 \dot z\notin\bigcup_{\gamma<\delta}\dot{\mathbb R}^{\mathbb{Q}_\gamma},
 \qquad
 \rho_{\leq(\xi,\theta_\xi)}(m,x,\dot z)=\lambda.
\]
Take \(p\leq q\) in the dense suborder and put
\(p_\delta=\pi_\delta(p)\). Since
\[
 s\leq p_\delta\quad\Longrightarrow\quad s\not\perp_{\mathbb{Q}} p,
\]
\(p_\delta\) forces both assertions over \(\mathbb{Q}_\delta\).
Thus \(\mathbb{Q}_\delta\restriction p_\delta\), with initial segments
\(\mathbb{Q}_\gamma\restriction\pi_\gamma(p)\), \(\gamma\leq\delta\),
is as required, also when \(\delta\) is a limit.
\end{proof}

\section{The uniformizing extension}\label{sec:uniformizing-extension}\label{section-L-final-uniformization}

\begin{theorem}\label{thm:uniformizing-extension}
There is a coding iteration \(\mathbb{P}\) of length \(\omega_1\) over \(W\)
such that, for every \(W\)-generic \(G\subseteq \mathbb{P}\), each \(A_m\)
is uniformized in \(W[G]\) by
\[
 U_m(x,y)\quad\Longleftrightarrow\quad
 (x,y)\in A_m\ \land\ \neg\Phi(\langle x,y,m\rangle).
\]
\end{theorem}

\begin{proof}
We use Lemmas~\ref{lem:coding-iterations}--\ref{lem:coded-reals}
and the absoluteness facts of
Subsection~\ref{subsec:persistence}.
Let \(\mathbb{P}_i\) be the initial segment of length \(i\), and put
\(G_i=G\cap \mathbb{P}_i\). By Lemma~\ref{lem:localization-of-reals},
\[
 \omega_1^{W[G_i]}=\omega_1,\qquad
 \mathbb R^{W[G_i]}=\{x_{i,j}:j<\omega_1\},\qquad
 \mathbb R^{W[G]}=\bigcup_{i<\omega_1}\mathbb R^{W[G_i]},
\]
where each enumeration is specified by a \(\mathbb{P}_i\)-name.
If \(\xi_i\) is the hierarchy in use at stage \(i\), abbreviate
\[
 \begin{aligned}
 \mathcal K_i&=K_{\xi_i}(W[G_i]),&
 \rho_i&=\rho_{\leq(\xi_i,\theta_{\xi_i})}^{W[G_i]},\\
 \Lambda_i(m,x)&=\Lambda^{\xi_i}_{\theta_{\xi_i}}(W[G_i],m,x),&
 j_i(m,x)&=j^{\xi_i}_{\theta_{\xi_i}}(W[G_i],m,x).
 \end{aligned}
\]
Let \(F_\xi\) be the triples explicitly forbidden in hierarchy
\(\xi\). List the following instructions cofinally often, using
\((i,j)\) for \(x_{i,j}\); skip references to initial segments not
yet constructed. Define \(y_*(m,x)\) at most once.

\begin{enumerate}
\item For \(b=(m,x,y)\in W[G_i]\),
\[
 \rho_i(b)\uparrow\quad\Longrightarrow\quad
 \exists \mathbb{Q}\in\mathcal K_i\ \exists q\in \mathbb{Q}\
 q\Vdash_{\mathbb{Q}}^{W[G_i]}\neg A(b).
\]
Append \(\mathbb{Q}\restriction q\).

\item If \(\mathcal K_i\) permits coding \((m,x,y)\), append a
coding stage for \(\langle x,y,m\rangle\).

\item Suppose \(y_*(m,x)\uparrow\) and
\(\Lambda_i(m,x)\ne\varnothing\). Put
\(\lambda=\min\Lambda_i(m,x)\). If \(\lambda\) is stable, set
\[
 y_*(m,x)=j_i(m,x)=J_{W[G_i]}(B),
\]
where \(B\) is the set of minimizers. Otherwise,
Lemma~\ref{lem:first-segment-minimizer} supplies
\(\mathbb{Q}\in\mathcal K_i\), \(\delta=\operatorname{lh}(\mathbb{Q})>0\), and
\(\dot y\) such that
\[
 \mathbb{Q}\Vdash
 \dot y\notin\bigcup_{\alpha<\delta}\dot{\mathbb R}^{\mathbb{Q}_\alpha}
 \quad\land\quad
 \rho_{\leq(\xi_i,\theta_{\xi_i})}(m,x,\dot y)=\lambda.
\]
Append \(\mathbb{Q}\), put \(y_*(m,x)=\dot y^H\) for its generic \(H\),
and start a new hierarchy with
\[
 F_{\xi_i+1}=F_{\xi_i}\cup\{(m,x,y_*(m,x))\}.
\]
The restriction applies before the next factor.
\end{enumerate}

The instructions are carried out using the uniformly definable predicates
of Lemma~\ref{lem:generic-recursion} and the forcing theorem.
On a maximal antichain \(\mathcal A\subseteq\mathbb{P}_i\), decide
which instruction applies and choose its witnesses, including
\(\dot{\mathbb{Q}}_a\) and \(\delta_a<\omega_1\) with
\(a\Vdash\operatorname{lh}(\dot {\mathbb{Q}}_a)=\delta_a\).
At the \(k\)-th stage within a tail, mix its \(k\)-th factor names;
at limits take the prescribed limits of initial segments.
Thus the length of an instruction may depend on the generic.
Between instructions extend the earlier choice functions coherently
to all new nonempty sets of reals. At a limit \(\zeta\) of hierarchies, use
\[
 F_\zeta=\bigcup_{\xi<\zeta}F_\xi,\qquad
 \mathcal K^\zeta_0(M)=\bigcap_{\xi<\zeta}K_\xi(M).
\]
Each new hierarchy starts before the next factor; earlier hierarchies
continue to determine previously assigned ranks.
If no tail is appended, insert a stage with trivial almost disjoint
factor after the choice or the start of a new hierarchy. For the starting stages of the
instructions, therefore,
\[
 \begin{gathered}
 \gamma_0=0,\qquad
 \gamma_{\nu+1}=\gamma_\nu+\ell_\nu\quad(0<\ell_\nu<\omega_1),\qquad
 \gamma_\delta=\sup_{\nu<\delta}\gamma_\nu\quad(\delta\text{ limit}),\\
 \nu\leq\gamma_\nu<\omega_1,\qquad
 \sup_{\nu<\omega_1}\gamma_\nu=\omega_1.
 \end{gathered}
\]
The last line follows from regularity of \(\omega_1\). Together with Lemma~\ref{lem:localization-of-reals}, it ensures that each instruction is
considered after its parameters have appeared.

\paragraph{Verification.}
For \(i<j<\omega_1\), the construction gives
\[
 \begin{gathered}
 \mathbb{P}_j/G_i\in\mathcal K_i,\\
 (m,x,y)\in F_{\xi_j}\setminus F_{\xi_i}
 \quad\Longrightarrow\quad y\notin W[G_i].
 \end{gathered}
\]
Thus Lemmas~\ref{lem:products-ranks}
and~\ref{lem:new-restrictions} yield
\[
 \rho_j\restriction(\omega\times\mathbb R^2)^{W[G_i]}=\rho_i.
\]
For \(b\in(\omega\times\mathbb R^2)^{W[G_i]}\), persistence,
the two absoluteness facts and instruction~(1) give
\[
 \begin{aligned}
 \rho_i(b)\downarrow
 &\Longrightarrow\ \forall j\in[i,\omega_1)\quad W[G_j]\models A(b)\\
 &\Longrightarrow\ W[G]\models A(b),\\[2pt]
 \rho_i(b)\uparrow
 &\overset{(1)}{\Longrightarrow}\ \exists j\in[i,\omega_1)\quad
                         W[G_j]\models\neg A(b)\\
 &\Longrightarrow\ W[G]\models\neg A(b).
 \end{aligned}
\]
Hence
\[
 W[G]\models A(b)\quad\Longleftrightarrow\quad\rho_i(b)\downarrow.
\]
Using the trivial forcing and preservation of the least attainable rank,
\[
 \begin{aligned}
 A_{m,x}^{W[G]}\ne\varnothing
 &\Longrightarrow\ \exists i,y\quad
       x,y\in W[G_i]\ \land\ \rho_i(m,x,y)\downarrow\\
 &\Longrightarrow\ \exists i\quad\Lambda_i(m,x)\ne\varnothing\\
 &\overset{(3)}{\Longrightarrow}\ y_*(m,x)\downarrow\ \land\
                    W[G]\models A(m,x,y_*(m,x)).
 \end{aligned}
\]

Fix \(A_{m,x}^{W[G]}\ne\varnothing\), and let
\(y_*=y_*(m,x)\) be chosen at \(\gamma\). Put
\[
 \lambda=\rho_\gamma(m,x,y_*),\qquad
 \mathcal C=\{u_i:i<\omega_1,\ i\text{ is a coding stage}\}.
\]
If \(y_*\) was chosen at a stable minimum, preservation gives,
for every \(i<\omega_1\),
\[
 \begin{aligned}
 x,y_*\in W[G_i]
 &\Longrightarrow\ \rho_i(m,x,y_*)\downarrow\\
 &\Longrightarrow\ \Lambda_i(m,x)\ne\varnothing\\
 &\Longrightarrow\ j_i(m,x)=j_\gamma(m,x)=y_*.
 \end{aligned}
\]
Therefore \eqref{eq:requirement} gives
\(\langle x,y_*,m\rangle\notin\mathcal C\).
Otherwise, instruction~(3) gives, also for limit \(\gamma\),
\[
 \begin{gathered}
 i<\gamma\ \Longrightarrow\ y_*\notin W[G_i],\qquad
 i\geq\gamma\ \Longrightarrow\ (m,x,y_*)\in F_{\xi_i},\\
 \Longrightarrow\quad \langle x,y_*,m\rangle\notin\mathcal C.
 \end{gathered}
\]

Let instruction~(2) act on \((m,x,y)\) at stage \(i\geq\gamma\),
where \(y\ne y_*\). In the stable case, the observation after
Lemma~\ref{lem:products-ranks} gives, for every level \((\eta,\beta)\) whose requirement \eqref{eq:requirement} applies at stage \(i\),
\[
 j^\eta_\beta(W[G_i],m,x)\downarrow
 \quad\Longrightarrow\quad j^\eta_\beta(W[G_i],m,x)=y_*.
\]
Otherwise, if \(\gamma_\nu\) starts the instruction choosing \(y_*\),
then \(y_*\notin W[G_{\gamma_\nu}]\). By
Lemma~\ref{lem:new-restrictions} over \(W[G_{\gamma_\nu}]\)
and the same observation,
\[
 \begin{gathered}
 \exists \mathbb{Q}\in\mathcal K_i\ \exists q\in \mathbb{Q}\ \exists\dot z\quad\\
 q\Vdash_{\mathbb{Q}}^{W[G_i]}
 \bigl(\dot z\notin\mathbb R^{W[G_i]}\ \land\
       \rho_{\leq(\xi_i,\theta_{\xi_i})}(m,x,\dot z)=\lambda\bigr),\\
 \Longrightarrow\quad j^\eta_\beta(W[G_i],m,x)\uparrow
 \quad\text{for every such level }(\eta,\beta).
 \end{gathered}
\]
Since
\[
 \{z:(m,x,z)\in F_{\xi_i}\}\subseteq\{y_*\},
\]
instruction~(2) yields
\[
 y\ne y_*\quad\Longrightarrow\quad
 u_i=\langle x,y,m\rangle\in\mathcal C.
\]
By Lemma~\ref{lem:coded-reals}, in \(W[G]\),
\[
 \Phi(\langle x,y,m\rangle)
 \quad\Longleftrightarrow\quad
 \langle x,y,m\rangle\in\mathcal C
 \quad\Longleftrightarrow\quad y\ne y_*.
\]
Together with \(U_m\subseteq A_m\), this gives
\[
 \{y:U_m(x,y)\}=
 \begin{cases}
  \{y_*(m,x)\},&A_{m,x}\ne\varnothing,\\
  \varnothing,&A_{m,x}=\varnothing.
 \end{cases}
\]
Finally,
\[
 A_m\in\Pi^1_k,\quad\Phi\in\Sigma^1_k
 \quad\Longrightarrow\quad U_m\in\Pi^1_k.
\]
Include real parameters in the first coordinate of a universal relation.
\end{proof}

\paragraph{Remark.}
The arguments of
Lemmas~\ref{lem:choice-functions}--\ref{lem:first-segment-minimizer}
and the construction above use only the iteration and coding properties
in Lemmas~\ref{lem:coding-iterations}--\ref{lem:coded-reals}
and the two absoluteness facts of
Subsection~\ref{subsec:persistence}. Thus the same argument
applies over another ground model once the corresponding properties
are verified. We use this in
Subsection~\ref{subsec:mn-uniformization-transfer}.

\begin{corollary}\label{cor:uniformization-over-L}
There is a cardinal-preserving set-generic extension of $L$ in which
$\CH$ and the $\Pi^1_3$-uniformization property hold.
\end{corollary}
\begin{proof}
First force with the branch product producing $W$, then with the
iteration of Theorem~\ref{thm:uniformizing-extension}. Cardinal preservation and
\(\CH\) follow from Theorem~\ref{properties-of-W} and
Lemma~\ref{lem:localization-of-reals}. Since \(\Phi\) is
\(\Sigma^1_3\), the uniformizing graphs are \(\Pi^1_3\).
\end{proof}

\section[Forcing over the canonical inner models Mn]{Forcing over the canonical inner models $M_n$}
\label{section-forcing-over-Mn}

Fix for the remainder of this section a natural number $n$ with
$1\leq n<\omega$, and assume that  $M_n^{\#}$ exists. As always we let $M_n$ be the resulting inner model one obtains when iterating the sharp out of the universe.
Throughout this section, $\omega_1$ denotes $\omega_1^{M_n}$.
Our goal is to prove the following theorem.

\begin{theorem*}
There is a set-generic extension of $M_n$ in which the
$\Pi^1_{n+3}$-uniformization property holds.
\end{theorem*}

The proof follows the construction over $L$ from the preceding sections, but
three ingredients require $M_n$-analogues: the recovery of the canonical lower
part $M_n|\omega_1$, the definition of the coding predicate, and
the absoluteness argument for witnesses to nonmembership in the sections.  We begin with the inner-model-theoretic background needed
for these analogues.

\subsection{Inner-model-theoretic background and generic absoluteness}
\label{subsection-Mn-background-and-absoluteness}

We use the notation and comparison conventions of Steel's analysis of mice
with finitely many Woodin cardinals; see \cite{Steel2,Steel3}.  An
\emph{ordinary premouse} means a premouse in the parameter-free premouse
language.  In particular, no real or other base set is included as a predicate.

\begin{definition}
\label{definition-n-smallness}
Let $\mathcal M$ be a premouse and let
$\eta<\operatorname{Ord}^{\mathcal M}$.  We say that $\mathcal M$ is
\emph{$n$-small above $\eta$} if, whenever $E$ is an extender on the
$\mathcal M$-sequence and
\[
   \eta<\operatorname{crit}(E),
\]
the initial segment
$\mathcal J^{\mathcal M}_{\operatorname{crit}(E)}$ does not have $n$ Woodin
cardinals above $\eta$.  We say that $\mathcal M$ is \emph{$n$-small} if it
is $n$-small above $0$.
\end{definition}

We use Steel's $\Pi_n$-iterability predicates \cite[Definitions~1.4 and~1.6]{Steel2}; their definition depends on the parity of $n$.
For a premouse $\mathcal M$ and a cutpoint $\eta$, write
$\operatorname{It}_n(\mathcal M,\eta)$ for this predicate above
$\eta$. We will use the cutpoint $\omega$.

Steel's analysis gives the following complexity statement.

\begin{fact}
\label{fact-projective-complexity-Pin-iterability}
For each fixed $1\leq n<\omega$, the relation
\[
   ``x\text{ codes a countable premouse which is }
   \Pi_n\text{-iterable above }\omega."
\]
is $\Pi^1_{n+1}$ in the real code $x$.  Moreover, every countable realizable
$n$-small premouse is $\Pi_n$-iterable at the indicated cutpoint.
\end{fact}

We next specify the generic extensions in which the comparison lemma will be
used.  All cardinalities in the following definition are computed in $M_n$.

\begin{definition}
\label{definition-small-Mn-generic-extension}
A model $N$ is a \emph{small $M_n$-generic extension} if
\[
   N=M_n[G]
\]
for some cardinal-preserving forcing $\forceP\in M_n$ such that
$M_n\models|\forceP|\leq\aleph_1$ and some $M_n$-generic filter
$G\subseteq\forceP$.  A \emph{small two-step extension} is an extension
\[
   M_n[G]\subseteq M_n[G][H]
\]
obtained from a two-step forcing
$\forceP*\dot{\forceQ}\in M_n$ which is forced to be cardinal preserving and
of size at most $\aleph_1$.
\end{definition}

The preparatory branch forcing and every forcing iteration used later will be
shown to have size $\aleph_1$ and to preserve cardinals.  Hence all models
which occur in the proof are small $M_n$-generic extensions in the sense of
Definition~\ref{definition-small-Mn-generic-extension}.

The second background ingredient is the mouse-operator closure which yields
projective generic absoluteness.  Recall that a set is
\emph{self-wellordered}, abbreviated \emph{swo}, if it codes a wellorder of
itself.  For a self-wellordered set $X$, let $M_k^{\#}(X)$ denote the canonical
active $X$-mouse with $k$ Woodin cardinals.  Thus $M_0^{\#}(X)$ is the usual
sharp over $X$.

\begin{fact}
\label{fact-Mn-closure-under-lower-mouse-operator}
The model $M_n$ is closed under the operator
\[
   X\longmapsto M_{n-1}^{\#}(X)
\]
below its least Woodin cardinal.  In particular, for every self-wellordered
$X\in H(\omega_2)^{M_n}$, the mouse $M_{n-1}^{\#}(X)$ exists and is normally
$\omega_2$-iterable in $M_n$.
\end{fact}

This closure is preserved by the small forcing notions which occur in the
construction.

\begin{lemma}
\label{lemma-small-forcing-preserves-lower-mouse-operator}
Let $N=M_n[G]$ be a small $M_n$-generic extension.  Then
\[
\begin{aligned}
   N\models{}&``M_{n-1}^{\#}(X)\text{ exists and is normally }
      \omega_2\text{-iterable}\\
   &\text{for every swo }X\in H(\omega_2)."
\end{aligned}
\]
Moreover, whenever $X$ belongs to both $M_n$ and $N$, the two models compute
the same mouse $M_{n-1}^{\#}(X)$.
\end{lemma}

\begin{proof}
Put \(\kappa=\omega_1^{M_n}\) and \(F(X)=M_{n-1}^{\#}(X)\).
For every swo \(X\in H(\kappa^+)^{M_n}\),
Fact~\ref{fact-Mn-closure-under-lower-mouse-operator} gives
\(F(X)\downarrow\) and normal \(\kappa^+\)-iterability.
Since \(|\forceP|^{M_n}\leq\kappa\) and \(\forceP\) preserves
cardinals, \cite[Lemma~3.4]{Schlicht}, with index \(n-1\), gives
the same assertion for swo \(X\in H(\kappa^+)^N\), as well as \(X\in H(\kappa^+)^{M_n},\ X\text{ swo} \quad\Longrightarrow\quad F^N(X)=F^{M_n}(X).\) \end{proof}

We can now state the precise generic absoluteness principle used later.  It is
the two-step form which permits real parameters introduced by the first
forcing.

\begin{theorem}
\label{theorem-small-two-step-projective-absoluteness}
Let
\[
   M_n[G]\subseteq M_n[G][H]
\]
be a small two-step extension.  Then
\[
   M_n[G]\prec_{\Sigma^1_{n+2}} M_n[G][H].
\]
Equivalently, for every real parameter $a\in M_n[G]$ and every
$\Sigma^1_{n+2}$ formula $\vartheta(v)$,
\[
   M_n[G]\models\vartheta(a)
   \quad\Longleftrightarrow\quad
   M_n[G][H]\models\vartheta(a).
\]
The same equivalence holds for $\Pi^1_{n+2}$ formulas.
\end{theorem}

\begin{proof}
Put \(N_0=M_n[G]\) and \(N_1=N_0[H]\). By
Lemma~\ref{lemma-small-forcing-preserves-lower-mouse-operator},
\(N_0\) is closed under \(X\mapsto M_{n-1}^{\#}(X)\), with
normal \(\omega_2\)-iterability, for swo \(X\in H(\omega_2)^{N_0}\).
With \(\kappa=\omega_1\) and \(|\forceQ|^{N_0}\leq\kappa\),
\cite[Lemma~3.7]{Schlicht} at index \(n-1\) gives \(N_0\prec_{\Sigma^1_{(n-1)+3}}N_1, \qquad\Sigma^1_{(n-1)+3}=\Sigma^1_{n+2}.\) For \(\psi\in\Pi^1_{n+2}\), apply this to
\(\neg\psi\in\Sigma^1_{n+2}\).
\end{proof}

The point at which Theorem~\ref{theorem-small-two-step-projective-absoluteness}
enters the uniformization proof is recorded in the following consequence.

\begin{corollary}
\label{corollary-persistence-Sigma-n-plus-3-witnesses}
Let
\[
   M_n[G]\subseteq M_n[G][H]
\]
be a small two-step extension, let $a\in M_n[G]$ be a real, and let
$\varphi(v)$ be a $\Sigma^1_{n+3}$ formula.  If
\[
   M_n[G]\models\varphi(a),
\]
then
\[
   M_n[G][H]\models\varphi(a).
\]
Consequently, if $A\subseteq(2^\omega)^2$ is $\Pi^1_{n+3}$ and
\[
   M_n[G]\models (x,y)\notin A,
\]
then $(x,y)\notin A$ remains true in $M_n[G][H]$.
\end{corollary}

\begin{proof}
Put \(N_0=M_n[G]\), \(N_1=M_n[G][H]\), and write
\(\varphi(a)\equiv\exists z\,\psi(z,a)\), where
\(\psi\in\Pi^1_{n+2}\). By
Theorem~\ref{theorem-small-two-step-projective-absoluteness},
\[
 \begin{aligned}
 N_0\models\varphi(a)
 &\Longrightarrow\exists z\in\mathbb R^{N_0}\ (N_0\models\psi(z,a))\\
 &\Longrightarrow\exists z\in\mathbb R^{N_0}\ (N_1\models\psi(z,a))\\
 &\Longrightarrow N_1\models\varphi(a).
 \end{aligned}
\]
For the last assertion take \(\varphi(x,y)\equiv(x,y)\notin A\).
\end{proof}

\subsection{The definable lower part of $M_n$}

\begin{definition}\label{def:lower-part}
Let $\mathcal I_n$ be the set of real codes for well-founded
ordinary premice $P$ such that
\[
 \begin{gathered}
 P\text{ is }\omega\text{-sound},\qquad \rho_\omega(P)=\omega,\\
 P\text{ is }n\text{-small},\qquad \operatorname{It}_n(P,\omega).
 \end{gathered}
\]
Here $\omega$ is a cutpoint, and we use Steel's iterability
predicate from \cite[Definitions~1.4 and~1.6]{Steel2}. We also
write $P\in\mathcal I_n$ when a real code for $P$ belongs to this
set. For a model $N$, let $\mathcal C_n^N$ be the class of transitive premice which have a code in $N$ satisfying this definition in $N$.
\end{definition}

The definition of $\mathcal I_n$ is parameter-free
$\Pi^1_{n+1}$ and is invariant under a change of the enumeration of
the coded premouse. Indeed the iterability predicate is
$\Pi_n$ over $\mathrm{HC}$ by \cite[Lemmas~1.5 and~1.7]{Steel2},
giving the bound in real codes in
Fact~\ref{fact-projective-complexity-Pin-iterability}.
The remaining checks use the satisfaction relations of the coded
structures and the finite-degree fine-structural definitions. They
do not increase the bound, and well-foundedness is $\Pi^1_1$.

\begin{lemma}\label{lem:lower-part}
Every member of $\mathcal C_n^{M_n}$
is an initial segment of $M_n$ of height below $\omega_1$, and
their heights are cofinal in $\omega_1$. If $N$ is a small
$M_n$-generic extension, then
\[
 \mathcal C_n^N=\mathcal C_n^{M_n},
 \qquad
 \bigcup\mathcal C_n^N=M_n|\omega_1.
\]
\end{lemma}

\begin{proof}
Write \(M=M_n\). By \cite[Theorem~2.3]{Steel2} and
\cite[Corollaries~3.12--3.13]{Steel3}, there is a family
\(\mathcal R\) of countable realizable sound segments \(Q\triangleleft M\)
such that
\[
 \rho_\omega(Q)=\omega\quad(Q\in\mathcal R),\qquad
 \mathbb R^M=\bigcup_{Q\in\mathcal R}(\mathbb R\cap Q),\qquad
 \sup_{Q\in\mathcal R}\operatorname{Ord}(Q)=\omega_1.
\]
Indeed, a bound \(\beta<\omega_1\) on the heights would give
\(\mathbb R^M\subseteq M|\beta\) with \(\lvert M|\beta\rvert^M\leq\aleph_0\).
By \cite[Corollary~4.4]{Steel2}, \(M\) is
\(\Sigma_n^{\mathrm{HC}}\)-correct. Hence, for \(Q\in\mathcal R\),
\cite[Lemma~2.2(2)]{Steel2} gives
\[
 V\models\operatorname{It}_n(Q,\omega)
 \quad\Longleftrightarrow\quad
 M\models\operatorname{It}_n(Q,\omega),
 \qquad\mathcal R\subseteq\mathcal C_n^M.
\]
For \(P\in\mathcal C_n^M\), choose \(Q\in\mathcal R\) with
\(\operatorname{Ord}(P)<\operatorname{Ord}(Q)\). These are
\(\omega\)-mice agreeing below the cutpoint \(\omega\), so
\cite[Lemma~2.2(3)]{Steel2} yields \(P\trianglelefteq Q\ \lor\ Q\trianglelefteq P \quad\Longrightarrow\quad P\trianglelefteq Q\triangleleft M.\) Let \(N\) be a small extension. For premouse codes \(p,q\) with
coded structures \(P_p,P_q\), put
\[
 \operatorname{Comp}(p,q)\ \Longleftrightarrow\
 (\exists Q\trianglelefteq P_q)\ P_p\cong Q
 \ \lor\ (\exists Q\trianglelefteq P_p)\ P_q\cong Q.
\]
An initial-segment cut and an isomorphism have real codes; thus
\(\operatorname{Comp}\in\Sigma^1_1\). If \(q\in M\) codes a
member of \(\mathcal C_n^M\), then
\[
 \chi(q)\equiv\exists p\,
  (p\in\mathcal I_n\land\neg\operatorname{Comp}(p,q))
 \ \in\Sigma^1_{n+2}.
\]
The comparison above gives \(M\models\neg\chi(q)\), and
Theorem~\ref{theorem-small-two-step-projective-absoluteness},
with trivial first forcing, gives \(N\models\neg\chi(q)\).

For \(P\in\mathcal C_n^N\), preservation of \(\omega_1\) gives
\(\operatorname{Ord}(P)<\omega_1\). Choose
\(Q\in\mathcal R\) above this height and a code \(q\in M\) for \(Q\).
Then \(N\models\neg\chi(q)\) implies \(P\trianglelefteq Q\).
Consequently \(P\) has a code \(p_0\in M\), and isomorphism
invariance followed by absoluteness gives \(N\models p_0\in\mathcal I_n \quad\Longleftrightarrow\quad M\models p_0\in\mathcal I_n.\) Thus \(\mathcal C_n^N\subseteq\mathcal C_n^M\). The reverse
inclusion follows by the same absoluteness applied to each ground-model
code. Finally,
\[
 \mathcal C_n^N=\mathcal C_n^M,\qquad
 \bigcup\mathcal C_n^N=\bigcup\mathcal R=M|\omega_1.
\]
\end{proof}

\paragraph{Canonical trees.}
Write $M=M_n$ and use $<_M$ for its canonical construction order.
The order on any initial segment is its own construction order,
and every predecessor of an element of that segment belongs to it
\cite{Steel2,Steel3}. All the following choices are made in $M$.

\begin{lemma}\label{lem:Mn-diamond}
There is a canonical $\diamondsuit_{\omega_1}$ sequence
$\langle a_\lambda:\lambda<\omega_1\rangle$ over $M$ such that, for every $\lambda<\omega_1$, all initial segments of $M$ containing a computation of $\langle a_\gamma:\gamma<\lambda\rangle$ compute the same sequence.
\end{lemma}

\begin{proof}
Put \(a_0=a_{\eta+1}=\varnothing\). At a nonzero limit
\(\lambda<\omega_1\), define
\[
 \begin{aligned}
 (X_\lambda,C_\lambda)=\min_{<_M}\{(X,C):{}&X\subseteq\lambda,
     \ C\subseteq\lambda\text{ club},\\[-2pt]
     &\forall\gamma\in C\ (a_\gamma\ne X\cap\gamma)\},
 \qquad a_\lambda=X_\lambda.
 \end{aligned}
\]
The set of candidates is nonempty: choose
\(\gamma_j\nearrow\lambda\) with \(j<\gamma_j\), and take \(C=\{\gamma_j:j<\omega\},\qquad X=\{j<\omega:j\notin a_{\gamma_j}\}.\) Then \(j\in X\cap\gamma_j\Longleftrightarrow j\notin a_{\gamma_j}\).
The condition on \((X,C)\) refers only to the previous values and to \(\lambda,X,C\).
Since \(<_M\) restricts to the construction order of each initial
segment, any segment containing the preceding computation and its
least witness makes the same choice.

For diamond, use Jensen's condensation argument in \(M\)
\cite[Theorems~5.1 and~5.4(1)]{Steel3}. Suppose \((X,C)\) is the
\(<_M\)-least counterexample at \(\omega_1\). Choose an
\(\omega\)-sound \(Q\triangleleft M\) with
\[
 X,C,\langle a_\lambda:\lambda<\omega_1\rangle\in Q,
 \qquad\rho_\omega(Q)=\omega_1^Q=\omega_1.
\]
For this projectum calculation, new subsets of \(\omega_1\) appear
cofinally below \(M|\omega_1^+\), supplying sound levels with
\(\rho_\omega(Q)\leq\omega_1\); see \cite[Sections~2 and~5]{Steel3}.
If \(\rho_\omega(Q)<\omega_1\), its hull together with the standard
parameters is countable in \(M\), contrary to \(\omega_1\in Q\).

Take a countable elementary hull \(H\prec Q\) containing the
displayed data and the standard parameters through a degree witnessing
\(\rho_\omega(Q)=\omega_1\). For its inverse collapse \(\pi\),
\[
 \pi:\bar Q\cong H\prec Q,\qquad
 \delta=H\cap\omega_1=\operatorname{crit}(\pi)=\rho_\omega(\bar Q).
\]
The second alternative of the condensation lemma, which involves an extender, does not occur. For the
corresponding degree-zero ultrapower \(j\), its critical point
\(\nu<\delta\) is countable in \(Q\). A surjection
\(s:\omega\twoheadrightarrow\nu\) in \(Q\) would satisfy
\[
 j(s)(k)=j(s(k))=s(k)\quad(k<\omega),\qquad
 j(\nu)=\operatorname{ran}(j(s))=\operatorname{ran}(s)=\nu,
\]
a contradiction. Hence \(\bar Q\triangleleft M\).

In \(\bar Q\), the pair \((X\cap\delta,C\cap\delta)\) is least
for the recursion below \(\delta\). The condition is absolute, and the construction order of \(\bar Q\) is an initial segment of \(<_M\), so this pair
is least in \(M\) as well. Therefore \(a_\delta=X\cap\delta,\qquad \delta\in C,\qquad a_\delta\ne X\cap\delta,\) the required contradiction.
\end{proof}

Use this diamond in the simultaneous construction of
Lemma~\ref{canonical-L-Suslin-sequence}, with the $<_M$-least choice
at every step. It gives independent normal splitting Suslin trees
$\vec S^{\,n}=\langle S^n_\xi:\xi<\omega_1\rangle$ in $M$. In particular, the
finite product argument there seals each guessed maximal antichain
in the appropriate product. At each level use the $<_M$-least
enumeration by an initial segment of $\omega$ to label its nodes.
These labels are used in the coding below.

\subsection[The preparatory branch extension Wn]{The preparatory branch extension $W_n$}
\label{subsection-preparatory-branch-extension-Mn}

We now pass to the ground model over which the coding iteration will be
performed.  Let
\[
   \vec S^{\,n}
   =\langle S^n_\xi\mid\xi<\omega_1\rangle
\]
be the independent sequence just constructed.  All products in
this subsection are computed in $M_n$.

\begin{definition}
\label{definition-branch-products-over-Mn}
For $A\subseteq\omega_1$ with $A\in M_n$, put
\[
   \operatorname{Br}^n_A
   :=\prod_{\xi\in A}^{\mathrm{fs}}\mathbb{S}^n_\xi.
\]
We write
\[
   \operatorname{Br}_n
   :=\operatorname{Br}^n_{\omega_1}
   =\prod_{\xi<\omega_1}^{\mathrm{fs}}\mathbb{S}^n_\xi.
\]
Let $H\subseteq\operatorname{Br}_n$ be $M_n$-generic.  For
$\xi<\omega_1$, define
\[
   b_\xi
   :=\bigcup\{p(\xi):p\in H\text{ and }
                     \xi\in\operatorname{dom}(p)\}.
\]
For $A\subseteq\omega_1$ in $M_n$, let
\[
   H_A:=H\cap\operatorname{Br}^n_A.
\]
Finally, put
\[
   W_n:=M_n[H].
\]
\end{definition}

Thus $b_\xi$ is the branch added through $S^n_\xi$ by the $\xi$-th
coordinate.  As before, to \emph{destroy} a tree $S^n_\xi$ means to add a
cofinal generic branch through it.

\begin{lemma}
\label{lemma-branch-product-calculus-Mn}
Let $A,B\in M_n$ satisfy
\[
   A\subseteq B\subseteq\omega_1.
\]
Then the following hold in $M_n$.
\begin{enumerate}
   \item $\operatorname{Br}^n_A$ is c.c.c.
   \item $\operatorname{Br}^n_A$ is a complete subforcing of
         $\operatorname{Br}^n_B$, and
         \[
            \operatorname{Br}^n_B
            \cong
            \operatorname{Br}^n_A
            \times
            \operatorname{Br}^n_{B\setminus A}.
         \]
   \item If $\xi\notin A$, then
         \[
            \operatorname{Br}^n_A\Vdash
            ``S^n_\xi\text{ is a Suslin tree}."
         \]
\end{enumerate}
\end{lemma}

\begin{proof}
Independence gives
\(\prod_{\xi\in e}S^n_\xi\text{ Suslin}\) for
\(\varnothing\ne e\in[\omega_1]^{<\omega}\) in \(M_n\).
Apply Lemma~\ref{finite-support-branch-products} there for (1).
The coordinate map \(p\longmapsto(p\restriction A,p\restriction(B\setminus A))\) gives (2), including \(\operatorname{Br}^n_A\lessdot\operatorname{Br}^n_B\).
For \(\xi\notin A\),
\[
 \operatorname{Br}^n_A\times\mathbb{S}^n_\xi
 \cong\operatorname{Br}^n_{A\cup\{\xi\}}\text{ is c.c.c.}
 \quad\Longrightarrow\quad
 \operatorname{Br}^n_A\Vdash``S^n_\xi\text{ is Suslin},"
\]
by the product characterization of Suslin-tree preservation.
\end{proof}

\begin{theorem}
\label{theorem-properties-of-Wn}
The model $W_n=M_n[H]$ has the following properties.
\begin{enumerate}
   \item The forcing $\operatorname{Br}_n$ is c.c.c. and has size
         $\aleph_1$ in $M_n$.  Consequently, $W_n$ preserves every cardinal
         and cofinality of $M_n$, and
         \[
            W_n\models\CH.
         \]

   \item Every $S^n_\xi$ is destroyed in $W_n$, as witnessed by the branch
         $b_\xi$.

   \item If $A\subseteq\omega_1$ belongs to $M_n$, then $H_A$ is
         $\operatorname{Br}^n_A$-generic over $M_n$, and, for every
         $\xi<\omega_1$,
         \[
            M_n[H_A]\models
            ``S^n_\xi\text{ has a cofinal branch}"
            \quad\Longleftrightarrow\quad
            \xi\in A.
         \]

   \item Every real in $W_n$ belongs to $M_n[H_A]$ for some countable set
         $A\subseteq\omega_1$ with $A\in M_n$.

   \item More generally, every member of $H(\omega_1)^{W_n}$ belongs to
         $M_n[H_A]$ for some countable set
         $A\subseteq\omega_1$ with $A\in M_n$.
\end{enumerate}
\end{theorem}

\begin{proof}
Apply Theorem~\ref{properties-of-W} inside \(M_n\), with
\[
 (L,\vec S,\operatorname{Br},W)
 \quad\mapsto\quad(M_n,\vec S^{\,n},\operatorname{Br}_n,W_n).
\]
Lemma~\ref{lemma-branch-product-calculus-Mn} supplies the same
product and preservation hypotheses; the argument with nice names from the proof of Theorem~\ref{properties-of-W} gives (4) and (5).
\end{proof}

The last two clauses are the branch-coordinate localization needed by the
coding argument.  We record a form in which a prescribed countable set of
coordinates may simultaneously be retained.

\begin{corollary}
\label{corollary-countable-branch-supports-Mn}
Let $a\in H(\omega_1)^{W_n}$, and let
$C\in[\omega_1]^\omega\cap M_n$.  Then there is
$A\in[\omega_1]^\omega\cap M_n$ such that
\[
   C\subseteq A,
   \qquad
   a\in M_n[H_A],
\]
and
\[
   M_n[H_A]\models
   ``S^n_\xi\text{ is a Suslin tree}"
   \qquad(\xi\in\omega_1\setminus A).
\]
In particular, no cofinal branch through any $S^n_\xi$ with
$\xi\notin A$ belongs to $M_n[H_A]$.
\end{corollary}

\begin{proof}
By Theorem~\ref{theorem-properties-of-Wn}(5), choose
\(A_0\in[\omega_1]^\omega\cap M_n\) with \(a\in M_n[H_{A_0}]\).
Put \(A=A_0\cup C\in[\omega_1]^\omega\cap M_n\). Then \(\operatorname{Br}^n_{A_0}\lessdot\operatorname{Br}^n_A, \qquad a\in M_n[H_{A_0}]\subseteq M_n[H_A].\) For \(\xi\notin A\), Lemma~\ref{lemma-branch-product-calculus-Mn}(3)
gives \(M_n[H_A]\models``S^n_\xi\text{ is Suslin},"\)
hence no cofinal branch through \(S^n_\xi\) in this model.
\end{proof}

\subsection[The coding over Mn]{The coding over $M_n$}

Write $m_\infty^n=M_n|\omega_1$ and use its fixed ordinal code, with
all premouse predicates, in $L[m_\infty^n,r]$. This code is in $M_n$,
and its transitive collapse recovers the lower part.

\paragraph{The almost disjoint family.}
Fix a recursive injection $\#\colon 2^{<\omega}\to\omega$. In $M$
define, for $\xi<\omega_1$,
\[
 x_\xi=\min_{<_M}\bigl(2^\omega\setminus\{x_\zeta:\zeta<\xi\}\bigr),
 \qquad
 d_\xi=\{\#(x_\xi\restriction k):k<\omega\}.
\]
Here $<_M$ is the canonical construction order on $M$. Each step is
defined because $\xi$ is countable and $M\models\CH$. The set $d_\xi$
is infinite, and if $\zeta<\xi$ and $k$ is the first position at which
$x_\zeta$ and $x_\xi$ differ, then
\[
 d_\zeta\cap d_\xi=\{\#(x_\xi\restriction j):j\leq k\}.
\]
Thus $D^n=\langle d_\xi:\xi<\omega_1\rangle$ is almost disjoint.
Write $D=D^n$ throughout this subsection. Its
initial segments are computed in initial segments of $M$ below
$\omega_1$, as shown below.

For \(Y\subseteq\omega_1\), the conditions of \(\mathbb A^{\,n}(Y)\) are pairs
\((s,F)\in[\omega]^{<\omega}\times[Y]^{<\omega}\), ordered by
\[
 (t,G)\leq(s,F)\quad\Longleftrightarrow\quad
 s\subseteq t,\quad F\subseteq G,\quad
 (t\setminus s)\cap\bigcup_{\xi\in F}d_\xi=\varnothing .
\]
For each stem \(s\), the set \(\{(s,F):F\in[Y]^{<\omega}\}\)
is centered, since \((s,\textstyle\bigcup_{j<\ell}F_j)\leq(s,F_j)\quad(j<\ell).\) For a generic filter \(G\) let \(r=\bigcup\{s:(s,F)\in G\}\). Density gives
\[
 \begin{aligned}
 \xi\in Y&\Longrightarrow\exists(s,F)\in G\ (\xi\in F)
                 \Longrightarrow |r\cap d_\xi|<\omega,\\
 \xi\notin Y&\Longrightarrow
    |d_\xi\setminus\textstyle\bigcup_{\zeta\in F}d_\zeta|=\aleph_0
                 \Longrightarrow |r\cap d_\xi|=\aleph_0.
 \end{aligned}
\]
In the second line, arbitrarily many such points can be added to any
stem. Thus \(Y=\{\xi<\omega_1:|r\cap d_\xi|<\omega\}\).

\paragraph{The coding sets.}
Fix in $M$ a bijection $\varrho_n\colon ({}^{<\omega_1}2)^M\to\omega_1$,
using $\CH$ in $M$.
Let $\mathbb C_{M_n}=({}^{<\omega_1}2)^M$, ordered by end-extension.
This is a dense suborder of $\operatorname{Add}(\omega_1,1)^M$. For a $\mathbb C_{M_n}$-generic $g$ over the
current model set
\[
 \begin{aligned}
 h(g)&=\{\varrho_n(g\restriction\eta):\eta<\omega_1\},\\
 \operatorname{Sel}(h,u)&=\{\tau_u(\gamma,k):\gamma\in h,\ k<\omega\},
 \qquad \tau_u(\gamma,k)=\omega\gamma+2k+u(k).
 \end{aligned}
\]
We call $h(g)$ the \emph{coding set} of $g$. Its cardinality is
$\omega_1$. For distinct $g$ and $g'$ the intersection
$h(g)\cap h(g')$ is bounded in $\omega_1$, since it consists of the
images of their common initial segments.

To code \(u\in2^\omega\), let \(g\) be the generic subset of \(\omega_1\) added at a new Cohen coordinate and put \(h=h(g)\). The forcing at the second step is
\(\mathbb A^{\,n}(Y_n(u,g))\), where \(Y_n(u,g)\subseteq\omega_1\) is defined below.
The set \(Y_n(u,g)\) codes \(u\), \(h\), and the branches \(b_\xi\), \(\xi\in\operatorname{Sel}(h,u)\).

\paragraph{Codes and restriction to countable heights.}
Fix a finite fragment of \(\ZFCminus\) containing the instances of separation, collection and replacement used below, the existence of satisfaction relations for set-sized structures and of transitive collapses of wellorders, and
add the axiom that \(\omega_1\) exists. A \emph{decoding model} is
a countable transitive model of this theory. Let
\(\langle\cdot,\cdot\rangle_{\rm ord}\) be G\"odel's
pairing function on ordinals, which enumerates pairs by increasing
maximum and lexicographically for a fixed maximum; for triples iterate
it. If \(N\) is a decoding model and \(\alpha=\omega_1^N\), then pairing maps \(\alpha\times\alpha\) onto \(\alpha\). Indeed, the part of the enumeration up to any pair of ordinals below \(\alpha\) is countable in \(N\).
Conversely the rank of a pair in the enumeration is at least the
maximum of its entries. Similarly, \(j\cdot\zeta+i<\alpha\) if and only if \(\zeta<\alpha\), for \(i<j<\omega\).
Also \(\tau_u(\gamma,k)<\alpha\) for \(\gamma<\alpha\) and
\(k<\omega\), since \(\alpha=\omega_1^N\).

Present the nodes of each \(S^n_\xi\) by pairs \((\lambda,i)\), where
\(\lambda\) is the level and \(i\) is the index in the \(<_M\)-least
bijection from an initial segment of \(\omega\) onto that level.
The bijection is chosen as soon as the level has been constructed.
We use these relabelled trees throughout, with the tree orders
transported along the bijections. In the codes below a node is identified with
the ordinal pairing of \(\lambda,i\). Thus the
restriction of the indexed tree to level \(<\alpha\) is obtained by
keeping exactly the node indices below \(\alpha\). We regard each
\(b_\xi\) as the set of indices of the nodes of its downward closure.

\paragraph{Localization.}
For \(\gamma\in h\) and \(k<\omega\) use the index
\(\tau_u(\gamma,k)\) and put
\[
\begin{aligned}
 B&=\{\langle\gamma,k,\zeta\rangle_{\rm ord}:
           \gamma\in h,\ k<\omega,\ \zeta\in b_{\tau_u(\gamma,k)}\},\\
 X&=\{3\cdot k:k<\omega,\ u(k)=1\}
       \cup\{3\cdot\gamma+1:\gamma\in h\}
       \cup\{3\cdot z+2:z\in B\}.
\end{aligned}
\]
The odd part of \(Y_n(u,g)\) codes \(X\). For the even part, define
a continuous increasing sequence in \(\operatorname{acc}(h)\) by
\[
 \begin{aligned}
 c_0&=\min\{\zeta\in\operatorname{acc}(h):\omega<\zeta\},\\
 c_{\nu+1}&=\min\{\zeta\in\operatorname{acc}(h):c_\nu+\omega<\zeta\},
 &c_\lambda&=\sup_{\nu<\lambda}c_\nu
                  \quad(\lambda\text{ a nonzero limit}).
 \end{aligned}
\]
Fix in \(M\) a family \(\langle w_\eta:\omega\leq\eta<\omega_1\rangle\)
with \(w_\eta\) coding a wellorder of \(\omega\) of type \(\eta\).
Define
\[
 \begin{aligned}
 E&=w_{c_0}\cup\{c_\nu+k:\nu<\omega_1,\ k\in w_{c_{\nu+1}}\},\\
 Y_n(u,g)&=\{2\cdot\zeta:\zeta\in E\}
           \cup\{2\cdot\zeta+1:\zeta\in X\}.
 \end{aligned}
\]
Here the wellorder codes use the fixed recursive coding of relations
on \(\omega\). As in the construction over $L$, the even part localizes the code, and the odd part codes the branch data. The set \(Y_n(u,g)\) is defined from \(g\), \(u\), the selected branches and parameters in \(M\).

Suppose \(N\) is a decoding model,
\(\alpha=\omega_1^N\), and \(Y_n(u,g)\cap\alpha\in N\).
Put \(C=\{c_\nu:\nu<\omega_1\}\); continuity and strict increase
make \(C\) club. If \(\alpha\notin C\), the even part yields
\[
 \begin{array}{c|c|c}
 \text{case}&\text{block of the even part}&\text{decoded wellorder}\\ \hline
 \alpha<c_0&[0,\omega)&w_{c_0}\\
 c_\nu<\alpha<c_{\nu+1}&[c_\nu,c_\nu+\omega)&w_{c_{\nu+1}}
 \end{array}
\]
In the second case, \(c_\nu\in N\) is countable there, so
\(c_\nu+\omega<\alpha\). The decoded relation is well-founded in the ambient model, hence in \(N\). Its collapse in \(N\)
would give, for \(\eta=c_0\) or \(c_{\nu+1}\), \(f\in N,\quad f:(\omega,w_\eta)\cong\beta, \quad\beta<\alpha<\eta.\) The same \(f\) is an isomorphism in the ambient model, so \(\beta=\eta\), a contradiction. These cases exhaust
\(\omega_1\setminus C\). Thus \(\alpha\in C\subseteq\operatorname{acc}(h),\qquad \sup(h\cap\alpha)=\alpha.\) Using the ordinal pairing and the node labels, \(N\) decodes from the odd part
\[
 \bigl(u,h\cap\alpha,
   \langle b_{\tau_u(\gamma,k)}\restriction\alpha:
                  \gamma\in h\cap\alpha,\ k<\omega\rangle\bigr).
\]
Each restricted branch has exactly one node at every level
\(\lambda<\alpha\).

\paragraph{Correct local sequences.}
The set $\mathcal I_n$ is given by Definition~\ref{def:lower-part}.
By Lemma~\ref{lem:lower-part}, the class $\mathcal C_n^N$ is the same for all small extensions $N$ used here.
For each \(\alpha<\omega_1\) put
\[
 D^\alpha=D\restriction\alpha,\qquad
 \vec S^{\,n,\alpha}=\langle S^n_\xi\restriction\alpha:\xi<\alpha\rangle,
\]
where the latter restriction is by levels. A triple
\((m,D',\vec S')\) consists of a nonempty well-founded premouse of
height \(\alpha=\operatorname{Ord}(m)\) which is the union of its
proper initial segments, together with sequences of the indicated
shape. The triple is \emph{correct} if for every \(\tau<\alpha\) there
is \(\beta\in(\tau,\alpha)\) such that \(m|\beta\in\mathcal I_n\) and
a computation in \(m|\beta\) agrees with \(D',\vec S'\) below
\(\tau\). This means agreement on the entries of \(D'\) of
index \(<\tau\) and on the trees of index \(<\tau\) restricted to
levels \(<\tau\), including node indices and orders. A
\emph{computation} is an element of \(m|\beta\) recording an initial
part of the constructions of \(D\) and \(\vec S^{\,n}\).

We show that correct triples agree with \(M\) and that they exist at club many heights.
For every initial segment \(M|\beta\),
\[
 <_{M|\beta}=<_M\restriction(M|\beta)^2,\qquad
 x\in M|\beta\ \Longrightarrow\ \{y:y<_Mx\}\subseteq M|\beta.
\]
A computation of \(D\) records \(x_\xi,d_\xi\) and checks
\[
 \begin{gathered}
 x_\xi\in2^\omega\setminus\{x_\zeta:\zeta<\xi\},\qquad
 d_\xi=\{\#(x_\xi\restriction j):j<\omega\},\\
 \forall y\in2^\omega\,
 (y<_Mx_\xi\Rightarrow\exists\zeta<\xi\ (y=x_\zeta)).
 \end{gathered}
\]
The quantifier ranges over the segment containing the computation.
All earlier candidates belong to that segment, so induction on
\(\xi\) gives its agreement with \(M\). For the diamond recursion,
record the least-pair witnesses from Lemma~\ref{lem:Mn-diamond}; the conditions on them are absolute once the preceding values are fixed, and
the same induction applies.

A computation of the trees records the orders below the current level, the antichains guessed by the diamond sequence, the branches at limit levels, and the level
bijections. If \(a_\lambda\) codes a maximal antichain
\(A\) in the product on a finite set \(e\), restricted below
\(\lambda\), the chosen branches must satisfy
\[
 \forall f\in{}^e\omega\ \exists t\in A\ \forall\xi\in e\
                  (t(\xi)\in b_{\xi,f(\xi)}).
\]
The branches for each tree must also cover all its preceding nodes.
These requirements refer only to the countable structures already
constructed, and witnesses exist by
Lemma~\ref{canonical-L-Suslin-sequence}. Hence a computation and
its least witness in \(M|\beta\) give the same choice in \(M\): \(\min_{<_{M|\beta}}\{b\in M|\beta:T(b)\} =\min_{<_M}\{b\in M:T(b)\},\) where \(T\) is the condition on the choice, with the preceding data fixed.

Each computation up to a countable stage is hereditarily countable in \(M\), so
\cite[Section~2]{Steel3} and Lemma~\ref{lem:lower-part} give \(\text{computation up to }\tau\ \in H_{\omega_1}^M =M|\omega_1=\bigcup\mathcal C_n^{M_n}.\) Choose \(f\in M\), \(f:\omega_1\to\omega_1\), such that \(\tau<f(\tau),\qquad M|f(\tau)\in\mathcal I_n,\) with the computations up to \(\tau\) belonging to \(M|f(\tau)\).
Well-foundedness is \(\Pi^1_1\); the computation checks are
arithmetic, and all quantifiers over segments and computations
range over countable codes. Thus the predicate \(\operatorname{Correct}_n\) expressing correctness is \(\Pi^1_{n+1}\), and
\[
 \operatorname{Correct}_n(m,D',\vec S')
 \ \Longrightarrow\
 (m,D',\vec S')=(M|\alpha,D^\alpha,\vec S^{\,n,\alpha}).
\]

In the current small extension choose \(H_\Theta\) sufficiently
large and a countable \(A\prec H_\Theta\) containing the parameters
in use, \(m_\infty^n,\omega_1,D,\vec S^{\,n},f\), the map
\(\lambda\mapsto M|\lambda\), and enumerations in \(M\) of each
\(M|\lambda\), \(0<\lambda<\omega_1\). For the collapse
\(\pi:A\cong N\) and \(\alpha=A\cap\omega_1\), these enumerations give
\[
 \begin{gathered}
 A\cap m_\infty^n=M|\alpha,\qquad \omega_1^N=\alpha,\qquad
 f``\alpha\subseteq\alpha,\\
 \pi(m_\infty^n)=M|\alpha,\qquad
 \pi(D)=D^\alpha,\qquad
 \pi(\vec S^{\,n})=\vec S^{\,n,\alpha}.
 \end{gathered}
\]
For \(\xi<\alpha\), the node indices of \(S^n_\xi\) in \(A\)
are exactly those below \(\alpha\); the collapse preserves them
and their order. Closure under \(f\) supplies the computations
required for correctness. A continuous cofinal chain of these
submodels gives a club of such heights \(\alpha\).

\begin{lemma}\label{lem:Mn-decoder}
There are parameter-free formulas
\[
 \Psi_n(r,u)\in\Pi^1_{n+2},
 \qquad
 \Phi_n(u)\ \longleftrightarrow\ \exists r\,\Psi_n(r,u),
 \qquad
 \Phi_n\in\Sigma^1_{n+3},
\]
with the following properties.
\begin{enumerate}
\item If \(r\) is the real added at a coding stage for \(u\), then
\(\Psi_n(r,u)\) holds in the resulting extension and in every
further small extension.
\item If \(\Psi_n(r,u)\) holds in a small extension, then
\(u\in L[m_\infty^n,r]\), and there is an unbounded \(h\subseteq\omega_1\)
in \(L[m_\infty^n,r]\) such that \(L[m_\infty^n,r]\) contains a cofinal
branch through \(S^n_\xi\) for every \(\xi\in\operatorname{Sel}(h,u)\).
\end{enumerate}
\end{lemma}

\begin{proof}
Let \(\operatorname{Dec}_n(c;r,u)\) assert that \(c\) codes a
decoding model \(N\) with \(r,u,(m,D',\vec S')\in N\),
\(\operatorname{Ord}(m)=\omega_1^N=\alpha\), and
\((m,D',\vec S')\) correct in the ambient model.
In \(N\), decode \((u',h',B')\) from the odd part of
\[
 Y_r^N=\{\xi<\alpha:|r\cap d'_\xi|<\omega\},\qquad
 B'_{\gamma,k}=\{\zeta:\langle\gamma,k,\zeta\rangle_{\rm ord}\in B'\}.
\]
Write \(\operatorname{Br}_\alpha(b,S)\) for the assertion that
\(b\subseteq S\) is linearly ordered by \(<_S\) and
\(\forall\lambda<\alpha\ (|b\cap S_\lambda|=1)\).
Define
\[
 \begin{aligned}
 \operatorname{Test}(c;r,u)\ \Longleftrightarrow\quad
 &u'=u\ \land\ \sup h'=\alpha\\
 &{}\land\ \forall\gamma\in h'\ \forall k<\omega\
 \operatorname{Br}_\alpha(B'_{\gamma,k},S'_{\tau_u(\gamma,k)}).
 \end{aligned}
\]
Include satisfaction relations in \(c\), checked by the arithmetic
Tarski clauses. Since \(\operatorname{Dec}_n\in\Pi^1_{n+1}\)
and \(\operatorname{Test}\) is arithmetic,
\[
 \Psi_n(r,u)\equiv\forall c\,
 (\operatorname{Dec}_n(c;r,u)\Rightarrow\operatorname{Test}(c;r,u))
 \ \in\Pi^1_{n+2}.
\]
Thus \(\Phi_n(u)\equiv\exists r\,\Psi_n(r,u)\) is
parameter-free \(\Sigma^1_{n+3}\).

Suppose \(r\) is added by \(\mathbb A^{\,n}(Y_n(u,g))\), possibly restricted below a condition. Put
\(Y_r=\{\xi<\omega_1:|r\cap d_\xi|<\omega\}\).
For every \(c\) satisfying \(\operatorname{Dec}_n(c;r,u)\),
correctness and localization yield
\[
 \begin{gathered}
 Y_r=Y_n(u,g),\qquad Y_r^N=Y_r\cap\alpha,\qquad
 \alpha\in\operatorname{acc}(h),\\
 (u',h',B')=(u,h\cap\alpha,B\cap\alpha),\qquad
 B'_{\gamma,k}=b_{\tau_u(\gamma,k)}\restriction\alpha
 \quad(\gamma\in h\cap\alpha,\ k<\omega).
 \end{gathered}
\]
Hence \(\operatorname{Test}(c;r,u)\) and therefore \(\Psi_n(r,u)\).
For any further small extension \(N_0\subseteq N_1\),
Theorem~\ref{theorem-small-two-step-projective-absoluteness} gives
\(N_0\models\Psi_n(r,u)\Longleftrightarrow N_1\models\Psi_n(r,u)\).

Conversely, assume \(\Psi_n(r,u)\). In \(L[m_\infty^n,r]\), the
coherent computations in \(m_\infty^n\), with bounds cofinal in
\(\omega_1\), recover \(D,\vec S^{\,n}\). Decode
\((\hat u,\hat h,\hat B)\) from \(Y_r\). Let \(\mathcal D_{r,u}\)
be the heights \(\omega_1^N\) of models coded by
\(c\) with \(\operatorname{Dec}_n(c;r,u)\); this set contains a club by the argument above for correct triples. At every \(\alpha\in\mathcal D_{r,u}\),
\[
 Y_r^N=Y_r\cap\alpha,\qquad
 (u',h',B')=(\hat u,\hat h\cap\alpha,\hat B\cap\alpha).
\]
The formulas \(\operatorname{Test}(c;r,u)\) therefore give \(\hat u=u\) and \(\sup\hat h=\omega_1\).
Fix \(\gamma\in\hat h\), \(k<\omega\), and
\(\xi=\tau_u(\gamma,k)\). Then
\[
 \alpha\in\mathcal D_{r,u},\ \alpha>\max\{\gamma,\xi\}
 \quad\Longrightarrow\quad
 \operatorname{Br}_\alpha(\hat B_{\gamma,k}\cap\alpha,
                          S^n_\xi\restriction\alpha).
\]
For any two nodes and any level there is such an \(\alpha\) above them, whence
\(\operatorname{Br}_{\omega_1}(\hat B_{\gamma,k},S^n_\xi)\), with
\(\hat B_{\gamma,k}\in L[m_\infty^n,r]\).
\end{proof}

\subsubsection*{Coding iterations and exactness}

For a set $E\in M_n$ of at most
$\aleph_1$ Cohen indices, put
\[
 \mathbb{C}_E^n=\operatorname{Fn}_{<\omega_1}(E\times\omega_1,2)^{M_n}.
\]
Use the same sets $E_i$ of Cohen indices as in the construction over
$L$, with the reindexing maps now taken in $M_n$. For $\beta\leq\omega_1$,
the forcing through the stages $i<\beta$ is
\[
 \mathbb{P}^n_\beta=\mathbb{C}^n_{E_{<\beta}}*\dot {\mathbb{R}}^n_\beta,
 \qquad E_{<\beta}=\bigcup_{i<\beta}E_i,
 \qquad
 \mathbb{R}^n_\beta=\mathop{*}^{\mathrm{fs}}_{i<\beta}\dot {\mathbb{A}}_i.
\]
At stage $i$, names over $\mathbb{P}_i^n$ specify whether to code a real, the
real $u_i$ to be coded, the selected coordinate in $E_i$, and the
ground-model isomorphism onto a restriction of $\mathbb{C}_{E_i}^n$.
If the stage codes a real, let $g_i$ be the Cohen subset of $\omega_1$
obtained at the selected coordinate after applying this isomorphism.
Then $\dot {\mathbb{A}}_i$ is a name for $\mathbb A^{\,n}(Y_n(u_i,g_i))$,
possibly restricted below a condition named over $\mathbb{P}_i^n*\mathbb{C}_{E_i}^n$.
If there is no coding at stage $i$, its second factor is trivial.

The proof of Lemma~\ref{lem:coding-iterations} applies over \(M_n\).
For disjoint \(E,F\in M_n\),
\[
 \mathbb{C}^n_{E\cup F}\cong \mathbb{C}^n_E\times \mathbb{C}^n_F,\qquad
 \pi_{\lambda\beta}(c,s)=
 \bigl(c\restriction(E_{<\beta}\times\omega_1),
                         s\restriction\beta\bigr).
\]
Given \(u_i,g_i\), the lower part and selected branches determine
\(Y_i=Y_n(u_i,g_i)\); the order of the next almost disjoint factor is determined by \(D^n\), \(Y_i\) and the condition below which it is restricted.
Thus the projections give \(\mathbb{P}^n_\beta\lessdot \mathbb{P}^n_\lambda\).

For restrictions, use the ground-model isomorphisms
\(\psi_c:\mathbb{C}_E^n\cong \mathbb{C}_E^n\restriction c\). They change each
coordinate only on a bounded set. Compose them with the stage
isomorphisms; the translated data then satisfy \((g_i',u_i',Y_i')=(g_i,u_i,Y_i).\) This yields the presentations below conditions and their commuting
initial-segment projections. A tail defined below a Boolean value
\(b\) is extended by mixing the names of its almost disjoint factors as
\(\dot {\mathbb{A}}_i^*=\dot {\mathbb{A}}_i\) on \(b\) and
\(\dot {\mathbb{A}}_i^*=\mathbf1\) on \(\neg b\).

For two tails, reindex their Cohen coordinates disjointly in \(M_n\)
and translate each factor's names along its initial-segment embeddings.
The data \(g_i,u_i,Y_i\) agree in either order, giving \(\mathbb{R}\times \mathbb{Q}\cong \mathbb{R}*\dot {\mathbb{Q}}\cong \mathbb{Q}*\dot {\mathbb{R}}.\) For a prescribed \(\mathbb{R}\times \mathbb{Q}_j\) and a continuation \(\mathbb{T}\) of
\(\mathbb{Q}_j\) below \(b\), retain both index maps and extend only
the \(\mathbb{Q}_j\)-map to the new coordinates. Use the image of \(b\) under \(\mathbb{Q}_j\lessdot \mathbb{R}\times \mathbb{Q}_j\) for the almost disjoint factors of \(\mathbb{T}\).
Then \(\mathbb{R}\times(\mathbb{Q}_j*\mathbb{T})\cong(\mathbb{R}\times \mathbb{Q}_j)*\dot {\mathbb{T}}\) over the prescribed product. The embeddings remain unchanged at
successors and, by finite support, at limits.

For \(\lambda\leq\omega_1\), put \(\mathbb{C}=\mathbb{C}^n_{E_{<\lambda}}\).
In \(M_n\), \(\mathbb{C}\) is countably closed and \(|\mathbb{C}|\leq\aleph_1\).
We may force with \(\mathbb{C}\) before \(\operatorname{Br}_n\). The countable-closure and
finite-product argument from the construction over \(L\) gives
\[
 M_n[K]\models\operatorname{Br}_n\text{ is c.c.c.},\qquad
 (\operatorname{Ord}^{\omega})^{W_n[K]}
          =(\operatorname{Ord}^{\omega})^{W_n}
\]
for generic \(K\subseteq \mathbb{C}\), the equality by Easton's lemma.
Each almost disjoint factor, including a restriction below a
condition, is Knaster. Hence \(\mathbb{R}^n_\lambda\) is c.c.c.\ and has
a dense set of finite functions with almost disjoint conditions
as values, of size \(\leq\aleph_1\). Consequently
\[
 \operatorname{dens}(\operatorname{Br}_n*\dot {\mathbb{P}}^n_\lambda)
       \leq\aleph_1,\qquad
 \operatorname{Br}_n*\dot {\mathbb{P}}^n_\lambda
       \text{ preserves cardinals and forces }\mathsf{CH}.
\]
The same calculation applies to products and to tails extended by trivial factors outside a Boolean value. For an initial-segment quotient, take dense suborders
and mixed names for an enumeration of its dense set. This gives
the size bound for the combined two-step forcing; thus each such
inclusion is a small two-step extension of \(M_n\).

For the localization of reals, take a nice \(\mathbb{R}^n_{\omega_1}\)-name
\(\dot z\) after forcing with the full Cohen product. Its
countable ordinal code lies in \(W_n\) by Easton's lemma, and
\[
 \operatorname{supp}(\dot z)=
 \bigcup\{\operatorname{dom}(p):(\check k,p)\in\dot z\}
 \subseteq\beta\quad\text{for some }\beta<\omega_1.
\]
The factors below \(\beta\) use only their own initial segments
and \(E_{<\beta}\), so the dense suborders of \(\mathbb{R}^n_\beta\) agree
after \(\mathbb{C}^n_{E_{<\beta}}\) and after the full Cohen product.
Evaluating the same name, as in Lemma~\ref{lem:localization-of-reals},
therefore gives
\[
 (2^\omega)^{W_n[G_{\omega_1}]}=
 \bigcup_{\beta<\omega_1}(2^\omega)^{W_n[G_\beta]},\qquad
 G_\beta=G_{\omega_1}\cap \mathbb{P}^n_\beta.
\]

\begin{theorem}\label{thm:Mn-coded-reals}\label{theorem-no-unwanted-Mn-codes}
In a generic extension of $W_n$ by a coding iteration of length at
most $\omega_1$,
\[
 \Phi_n(u)\quad\Longleftrightarrow\quad
 \exists i\ (i\text{ is a coding stage and }u_i=u).
\]
\end{theorem}

\begin{proof}
At a coding stage \(i\), the generic real \(r_i\) for
\(\mathbb A^{\,n}(Y_n(u_i,g_i))\), also below a condition, satisfies
\(\Psi_n(r_i,u_i)\), hence \(\Phi_n(u_i)\), by
Lemma~\ref{lem:Mn-decoder}. The selected \(g_i\) is obtained
from the specified ground-model isomorphism, and the same witness
persists to the final extension.

Conversely, suppose \(\Psi_n(r,u)\) and no stage codes \(u\).
Let \(K\) be generic for the full Cohen product, and work in
\(W_n[K]\). Choose nice names \(\dot r,\dot u\) using the dense
suborder of finite functions. Let \(D_i\subseteq i\) collect the supports of the names for the parameters of stage \(i\): whether to code and which real,
the selected Cohen coordinate and isomorphism, and the condition
restricting the almost disjoint factor. Put
\[
 \begin{aligned}
 J_0&=\operatorname{supp}(\dot r)\cup\operatorname{supp}(\dot u),\\
 J_{\ell+1}&=J_\ell\cup\bigcup_{i\in J_\ell}D_i,
 &J&=\bigcup_{\ell<\omega}J_\ell\in[\omega_1]^{\leq\omega}.
 \end{aligned}
\]
Since \(D_i\subseteq J\cap i\) for \(i\in J\), induction gives
\[
 \mathbb{R}^n_{J\cap i}\lessdot \mathbb{R}^n_i,\qquad
 \dot r^{G_J}=r,\quad\dot u^{G_J}=u,\qquad G_J=G\cap \mathbb{R}^n_J.
\]
At successors the name for the factor uses only earlier stages in \(J\);
at limits the embeddings commute with the finite-support projections.
The names have a common countable ordinal code \(z\). Easton's
lemma and the c.c.c.\ of the branch product give
\[
 z\in W_n,\qquad z\in M_n[H_B]
 \quad\text{for some }B\in[\omega_1]^{\leq\omega}\cap M_n.
\]

Let \(\mathcal C=\{i\in J:i\text{ is a coding stage}\}\) and
\(I_i=\operatorname{Sel}(h_i,u_i)\) for \(i\in\mathcal C\), with
\(I_i=\varnothing\) otherwise. Distinct Cohen coordinates differ on an unbounded set, and the isomorphisms change them only on bounded sets;
thus \(i\ne j\in\mathcal C\ \Longrightarrow\ \sup(h_i\cap h_j)<\omega_1.\) Choose \(\eta<\omega_1\) such that
\(B\subseteq\eta\) and \(h_i\cap h_j\subseteq\eta\) for all
distinct \(i,j\in\mathcal C\). Lemma~\ref{lem:Mn-decoder} supplies
an unbounded \(h_r\) and branches in \(L[r,m_\infty^n]\) through
every \(S^n_{\tau_u(\gamma,k)}\), \(\gamma\in h_r\), \(k<\omega\).
Choose \(\gamma\in h_r\setminus(\eta+1)\); then
\(|\{i\in\mathcal C:\gamma\in h_i\}|\leq1\).
If this set is \(\{i\}\), choose \(k\) with \(u_i(k)\ne u(k)\);
otherwise put \(k=0\). Consequently
\[
 \xi=\tau_u(\gamma,k)\notin B\cup\bigcup_{i\in J}I_i.
\]

Put \(M_\xi=M_n[K][H\restriction(\omega_1\setminus\{\xi\})]\);
then \(z\in M_\xi\). As in Lemma~\ref{lem:coded-reals}, interpret
the names on \(J\) over \(M_\xi\), replacing every factor from the
first stage requiring the branch at \(\xi\) onward by the trivial
forcing. Denote the result by \(\widetilde {\mathbb{R}}_J^\xi\).
Before that stage, all data defining \(Y_n\) lie in the model.
Induction on the common dense suborders of finite functions gives,
in \(W_n[K]\),
\[
 \begin{gathered}
 d_\xi=\bigwedge_{i\in J}
       \|\xi\notin\dot I_i\|_{\operatorname{RO}(\mathbb{R}^n_J)},\\
 \widetilde {\mathbb{R}}_J^\xi\lessdot \mathbb{R}^n_J,\qquad
 \widetilde {\mathbb{R}}_J^\xi\restriction d_\xi
       \simeq \mathbb{R}^n_J\restriction d_\xi.
 \end{gathered}
\]
At successors, membership and order in the next factor agree in
common generic extensions of the preceding stages. At limits, finite
supports are bounded and the corresponding Boolean values are meets.
In the given extension, \(\xi\notin I_i\) for all \(i\in J\);
the induced \(M_\xi\)-generic \(\widetilde G\) therefore satisfies
\(r,u\in M_\xi[\widetilde G]\).

The Cohen product can precede the branch forcing, so
\(M_\xi\models``S^n_\xi\text{ is Suslin}."\)
The Knaster factors and finite support give
\[
 M_\xi\models\widetilde {\mathbb{R}}_J^\xi\times\mathbb{S}^n_\xi\text{ is c.c.c.}
 \quad\Longrightarrow\quad
 M_\xi[\widetilde G]\models``S^n_\xi\text{ is Suslin}."
\]
But \(m_\infty^n\in M_\xi\), whence
\(L[r,m_\infty^n]\subseteq M_\xi[\widetilde G]\); the branch
through \(S^n_\xi\) decoded there contradicts this conclusion.
\end{proof}

\subsection{Transfer of the uniformization construction}
\label{subsec:mn-uniformization-transfer}
Fix a recursive universal enumeration of $\Pi^1_{n+3}$ formulas
with two free real variables, and write $A^n_m$ for the corresponding
relations. We use the construction in the proof of
Theorem~\ref{thm:uniformizing-extension} with $W_n$, $A_m^n$, and
$\Phi_n$ in place of $W$, $A_m$, and $\Phi$. We verify the required
properties below.

The preceding arguments give the complete embeddings, products,
restrictions, and continuations used in
Theorem~\ref{thm:uniformizing-extension}. The density bound supplies
a set of names for tails and choice functions. If \(N\) is the
common initial extension and \(H_{\mathbb{R}},H_{\mathbb{Q}}\) are the two quotient
generics, then \(N[H_{\mathbb{R}}]\cap N[H_{\mathbb{Q}}]=N.\) Thus the intersection lemma applies to both the choice functions
and the restrictions on coding reals outside \(N\). The classes of allowable iterations are defined from the forcing relation for \(A_m^n\) and these embeddings.

For persistence write
\[
 \neg A_m^n(x,y)\ \Longleftrightarrow\
 \exists z\,\theta_m(z,x,y),\qquad\theta_m\in\Pi^1_{n+2}.
\]
For the extension \(N_{\mathbb{R}}=N[H_{\mathbb{R}}]\) by one factor and the product extension
\(N_{\mathbb{R},\mathbb{Q}}=N[H_{\mathbb{R}},H_{\mathbb{Q}}]\),
Theorem~\ref{theorem-small-two-step-projective-absoluteness} gives
\[
 \begin{aligned}
 N_{\mathbb{R}}\models\theta_m(z,x,y)
 &\Longrightarrow N_{\mathbb{R},\mathbb{Q}}\models\theta_m(z,x,y)\\
 &\Longrightarrow N_{\mathbb{R},\mathbb{Q}}\models\neg A_m^n(x,y).
 \end{aligned}
\]
This also applies after appending a tail forcing nonmembership to
one factor, including after further coding restrictions.
Conversely, if a selected value \(y\) satisfies \((x,y)\notin A_m^n\) in the final model, localization
and the same absoluteness give
\[
 \begin{aligned}
 W_n[G_{\omega_1}]\models\neg A_m^n(x,y)
 &\Longrightarrow\exists\beta<\omega_1\ \exists z\in W_n[G_\beta]\\[-2pt]
 &\hspace{12mm}W_n[G_\beta]\models\theta_m(z,x,y),
 \end{aligned}
\]
with \(x,y\in W_n[G_\beta]\). All these inclusions are small
two-step extensions by the size calculation. Hence the arguments
for persistence, least attainable ranks, their minimizers, the first
segment containing a new minimizer, and preservation of the selected
value transfer to \(M_n\).

\begin{theorem}\label{theorem-force-Pi-n-plus-3-uniformization-over-Mn}
Let $1\leq n<\omega$ and suppose that $M_n^\#$ exists. There is a
cardinal-preserving set-generic extension of $M_n$ satisfying
$\mathsf{CH}$ and the $\Pi^1_{n+3}$-uniformization property.
\end{theorem}

\begin{proof}
Force first with \(\operatorname{Br}_n\), then apply
Theorem~\ref{thm:uniformizing-extension} over \(W_n\), using the
properties verified above. Let \(N=W_n[G_{\omega_1}]\) and
\(\mathcal C=\{u_i:i<\omega_1,\ i\text{ is a coding stage}\}\).
In \(N\), put
\[
 U_m^n(x,y)\ \Longleftrightarrow\
 A_m^n(x,y)\land\neg\Phi_n(\langle x,y,m\rangle).
\]
For \(A^n_{m,x}\ne\varnothing\), the construction selects
\(y_*(m,x)\in A^n_{m,x}\). By Theorem~\ref{thm:Mn-coded-reals},
for every \(y\in A^n_{m,x}\),
\[
 \neg\Phi_n(\langle x,y,m\rangle)
 \ \Longleftrightarrow\ \langle x,y,m\rangle\notin\mathcal C
 \ \Longleftrightarrow\ y=y_*(m,x).
\]
Thus \(\{y:U_m^n(x,y)\}=\{y_*(m,x)\}\) for a nonempty section,
and equals \(\varnothing\) otherwise. Since
\(A_m^n\in\Pi^1_{n+3}\) and \(\Phi_n\in\Sigma^1_{n+3}\) by
Lemma~\ref{lem:Mn-decoder}, \(U_m^n\in\Pi^1_{n+3}\).
Include real parameters in the first coordinate by recursive pairing.
Preservation gives \(\operatorname{Card}^N=\operatorname{Card}^{M_n}\)
and \(N\models\mathsf{CH}\).
\end{proof}

\section{Further possible applications and open problems}
We close with several directions suggested by the construction.  First, the method
should have an analogue in the generalized Baire-space context at a regular
uncountable cardinal $\kappa$.  In particular, one may ask whether the
$\Pi^1_1$-uniformization problem for ${}^\kappa\kappa$ can consistently have
a positive solution.

The methods in this paper leave the following question open:
\begin{question}
Can the $\Pi^1_n$-uniformization property be forced over $L$ for $n >3$?
\end{question}

The method introduced here is so far limited to local effects.  It would be
interesting to force a less local, or even global, behaviour:
\begin{question}
Given a pair $n,m\in\omega$ such that $n\neq m$ and $n\neq m+1$, can one force
a universe in which the $\Pi^1_n$- and the $\Pi^1_m$-uniformization properties
both hold?
\end{question}

Finally, it seems worthwhile to combine the present method with other forcing
constructions in order to obtain models with properties that are usually
consequences of determinacy assumptions.  We record one paradigmatic question:
\begin{question}
Assume some sufficiently big large cardinal notion which can exist in $L$.  Is there a model which satisfies the
$\Pi^1_3$-uniformization property and in which every projective set of reals is
Lebesgue measurable?
\end{question}

\end{document}